\documentclass[preprint]{elsarticle}
\usepackage[dvipsnames]{xcolor}
\usepackage{amsmath,amsthm,amsfonts,amssymb,tikz}
\usepackage[hidelinks]{hyperref}
\usepackage[capitalise,noabbrev]{cleveref}
\usepackage{enumitem}

\journal{arXiv}

\hypersetup{colorlinks,linkcolor={red!50!black},citecolor={blue!50!black},urlcolor={blue!80!black}}

\newtheorem{thm}{Theorem}\crefname{thm}{Theorem}{Theorems}
\newtheorem{lem}[thm]{Lemma}\crefname{lem}{Lemma}{Lemmas}
\newtheorem{ex}[thm]{Example}\crefname{ex}{Example}{Examples}
\newtheorem{cor}[thm]{Corollary}\crefname{cor}{Corollary}{Corollaries}
\crefname{conj}{Conjecture}{Conjectures}
\newtheorem{prop}[thm]{Proposition}\crefname{prop}{Proposition}{Propositions}
\newtheorem{defn}[thm]{Definition}\crefname{defn}{Definition}{Definitions}
\newtheorem{rem}[thm]{Remark}\crefname{rem}{Remark}{Remarks}
\newtheorem{que}[thm]{Question}\crefname{que}{Question}{Questions}
\crefname{figure}{Figure}{Figure}

\newcommand{\atm}[2]{at^{#1}(#2)}
\newcommand{\F}[3]{\Omega^{#3}_{#2}(#1)}
\DeclareMathOperator{\NE}{\textrm{NE}}
\DeclareMathOperator{\rp}{\textrm{rp}}
\DeclareMathOperator{\red}{\textrm{red}}
\DeclareMathOperator{\rk}{\textrm{rk}}
\DeclareMathOperator{\at}{\textrm{at}}

\numberwithin{equation}{section}
\numberwithin{figure}{section}
\numberwithin{thm}{section}

\begin{document}

\begin{frontmatter}
\title{On the {M}\"obius Function and Topology of General Pattern Posets}

\author{Jason P. Smith\fnref{fn1}}
\address{Department of Computer and Information Sciences, University of Strathclyde, \\
             Glasgow, UK.}
\ead{jason.p.smith@strath.ac.uk}
\fntext[fn1]{This research was supported by the EPSRC Grant EP/M027147/1.}

\begin{abstract}
We introduce a formal definition of a pattern poset which encompasses several previously studied posets in the literature. Using this definition we present some general results on the M\"obius function and topology of such pattern posets. We prove our results using a poset fibration based on the embeddings of the poset, where embeddings are representations of occurrences. We show that the M\"obius function of these posets is intrinsically linked to the number of embeddings, and in particular to so called normal embeddings. We present results on when topological properties such as Cohen-Macaulayness and shellability are preserved by this fibration. Furthermore, we apply these results to some pattern posets and derive alternative proofs of existing results, such as Bj\"orner's results on subword order.
\end{abstract}
\end{frontmatter}

\section{Introduction}\label{sec:intro}

Pattern occurrence, or more generally the presence of substructures, has been studied on a wide range of combinatorial objects with many different definitions of a pattern; see \cite{Kit11} for an overview of the field. In many of these cases we can use the notion of pattern containment to define a poset on these objects, for example the classical permutation poset. Whilst many such pattern posets have been studied in isolation, there is no general framework for the study of these posets. Yet many of the known results follow a similar theme. By introducing a formal definition of a pattern poset we develop tools for studying these posets, which leads to some general results that helps in understanding why different pattern posets often have a similar structure.

We say a word $\alpha$ occurs as a \emph{$\rho$-pattern} in a word $\beta$ if there is a subsequence of~$\beta$ satisfying certain conditions $\rho$ with respect to $\alpha$, and we  call this subsequence an \emph{occurrence of $\alpha$}. We can define a binary relation $\alpha\le_\rho\beta$ if $\alpha$ occurs as a $\rho$-pattern in~$\beta$. If $\le_\rho$ satisfies the partial order conditions, then we can define a \emph{pattern poset} on the set of words in question. A variety of different pattern posets have been studied in the literature, where the main focus is to answer questions on the structure and topology of a pattern poset $P$ and its \emph{intervals}, that is, the induced subposets $[\alpha,\beta]=\{\lambda\in P\,|\,\alpha\le\lambda\le~\beta\}$.
 
The topology of a poset is considered by mapping the poset to a simplicial complex, called the \emph{order complex}, whose faces are the chains of the poset, that is, the totally ordered subsets. We refer the reader to \cite{Wac07} for an overview of poset topology. A poset is \emph{shellable} if its maximal chains can be ordered in a certain way. Shellability implies a poset has many nice properties, such as Cohen-Macaulayness. We define a poset $P$ to be \emph{Cohen-Macaulay} if the order complex of $P$, and of every interval of $P$, is homotopically equivalent to a wedge of top-dimensional spheres.

The study of patterns in words has received a lot of attention. Perhaps the simplest type of pattern in a word is that of \emph{subword order}, that is, $u=u_1\ldots u_a$ occurs as a pattern in $w=w_1\ldots w_b$ if there is a subsequence $w_{i_1}\ldots w_{i_a}$ such that $w_{i_j}=u_j$ for all~$j=1,\ldots,a$. Bj\"orner \cite{Bjo90} presented a formula for the M\"obius function of the poset of words with subword order and showed that this poset is shellable. The poset of words with \emph{composition order} has the partial order $u\le w$ if there is a subsequence $w_{i_1}\ldots w_{i_a}$ such that $u_i\le w_{i_j}$ for all $j=1,\ldots,a$. A formula for the M\"obius function of this poset is given by Sagan and Vatter in~\cite{SagVat06}. Furthermore, the poset of generalised subword order is considered by Sagan and Vatter in~\cite{SagVat06} and McNamara and Sagan in~\cite{McnSag12}. There are many other examples of studies of patterns in words and word-representable objects, such as in permutations, set partitions, trees, mesh patterns and more.

Perhaps the most studied pattern posets in recent years is that of permutation patterns, where a \emph{permutation} is a word on the alphabet of nonnegative integers with no repeated letters. The classical, and most studied, definition of a pattern in a permutation says that~$\sigma$ occurs as a pattern in the permutation~$\pi$ if there is a subsequence of~$\pi$ whose letters have the same relative order of size as the letters of $\sigma$. For example, $213$ occurs as a pattern in $35142$ in the subsequence $314$. The permutation pattern poset has been studied extensively but a complete understanding has proved elusive due to its complex nature. Some formulas for the M\"obius function of certain classes of permutations and certain properties of the topology have been given in~\cite{SagVat06,SteTen10,BJJS11,Smith13,McSt13,Smith14,Smith15}.

Many of the known results on the M\"obius function of pattern posets, including those mentioned above, depend on the number of \emph{normal embeddings}, defined in various but similar ways.  For example, they play an important role in the study of many different classes of intervals of the classical permutation poset. We define an \emph{embedding} of $\alpha$ in $\beta$ as a sequence of dashes and the letters of $\alpha$, such that the positions of the non-dash letters give an occurrence of $\alpha$ in $\beta$ and deleting all the dashes results in $\alpha$.  The definition of when an embedding is normal varies, but all follow a similar theme. Perhaps the simplest definition is that of Bj\"orner's for subword order \cite{Bjo90}, where the normal condition is that the only positions that can be dashes are the leftmost positions in the maximal consecutive sequences of equal letters. For example, $1-2-1$ is a normal embedding of $121$ in $12211$.

We introduce a simple definition for normal embeddings which extracts the common theme from those in the literature. Using this definition we prove that the M\"obius function of a pattern poset, satisfying certain restrictions, equals the number of normal embeddings, plus an extra term that we describe explicitly. This extra term embodies the variations in the many definitions of normal embeddings. Intriguingly, this extra term often vanishes, or can be shown to be zero, which allows us to compute the M\"obius function in polynomial time. Furthermore, our general result can be used to prove many of the existing results on the M\"obius function of various pattern posets.

Poset fibrations are instrumental to our results. A fibration of a poset $Q$ consists of another poset~$P$, called the \emph{total space}, and a rank and order preserving surjective map~$f:P\rightarrow Q$. Poset fibrations were first studied by Quillen in \cite{Qui78} and a good overview is given in \cite{Bjo05}. It was shown by Quillen that Cohen-Macaulayness is maintained across a poset fibration satisfying certain conditions. We prove that in some case shellability can also be preserved in a similar manner.

We introduce a poset fibration on an interval $[\alpha,\beta]$ of a pattern poset. The total space of this fibration is built from the embeddings of $\lambda$ in $\beta$, for all~${\lambda\in(\alpha,\beta)}$. This total space has a much nicer structure than the original poset, which allows us to compute the M\"obius function and topology of the total space. We can then use known results on poset fibrations to get results for the original interval.

It is known that a poset is not shellable if it contains a disconnected subinterval of rank greater than $2$. We say an interval is \emph{zero split} if its set of embeddings can be partitioned into two parts $A$ and $B$ such that no position appears as a dash in an embedding from $A$ and in an embedding from $B$. In~\cite{McSt13} it is shown that an interval of the classical permutation poset is disconnected if and only if it satisfies a slightly stronger condition than being zero split. We introduce a definition of \emph{strongly zero split} which generalises this result to pattern posets. This implies that if an interval contains a strongly zero split subinterval of rank greater than~$2$, then it is not shellable. 

In \cref{sec:notation} we introduce some notation used throughout the paper and in \cref{sec:patPos} we introduce pattern posets. In \cref{sec:posfib} we introduce two poset fibrations on pattern posets. In \cref{sec:results} we apply these poset fibrations to prove some results on the  M\"obius function and topology of pattern posets. In \cref{sec:app} we apply these results to the poset of words with subword order, which provides an alternative proof of Bj\"orner's result on the M\"obius function of this poset. We also consider the consecutive permutation pattern poset and provide an alternative proof for the results on the M\"obius function given in \cite{BFS11} and \cite{SW12}. Finally, in \cref{sec:FW} we propose some question for future work.

\section{Notation and Preliminaries}\label{sec:notation}
We begin by introducing some necessary notation on words and posets. For further background on words see \cite{Kit11} and for further background on posets see \cite[Chapter 3]{Sta97}.
\begin{defn}
A \emph{word} is a sequence of letters from an alphabet $\Sigma$. The \emph{length} of a word $w$, denoted~$|w|$, is the number of letters in the word and we use $w_i$ to denote the letter in position $i$ of $w$. We denote the set of words on the alphabet~$\Sigma$ by $\Sigma^*$.
\end{defn}
\begin{ex}
If $\Sigma=\{0,1\}$ then $01001$ is a word of length $5$ on the alphabet~$\Sigma$.
\end{ex}

Note that we use the convention that the first position is number $1$, not $0$. We can apply many different restrictions to words to get different combinatorial objects. We are particularly interested in permutations, which can be defined in the following way:
\begin{defn}
 A \emph{permutation} is a word with no repeated letters on the alphabet of positive integers. The \emph{reduced form} of a permutation $\sigma$ is the permutation  $\red(\sigma)$, where if $\sigma_i$ is the $k$'th largest valued letter then $\red(\sigma)_i=k$. Two permutations are considered to be the same if they have the same reduced form.
\end{defn}
\begin{ex}
If $\sigma=264$, then $\red(\sigma)=132$.
\end{ex}

When studying the M\"obius function or topology of a poset it is often necessary that there is a unique minimal and unique maximal element, called the \emph{bottom} and \emph{top} elements, respectively.

\begin{defn}
A poset $P$ is \emph{bounded} if it has a unique bottom and top element, which we denote~$\hat{0}$ and $\hat{1}$, respectively. If $P$ is not bounded we create the bounded poset $\hat{P}$ by adding a bottom and top element. The interior of a bounded poset is obtained  by removing the bottom and top elements.
\end{defn}

Many of the posets that we look at are infinite, so it makes sense to limit our investigation to smaller subposets.

\begin{defn}
An \emph{interval} of a poset $P$ is an induced subposet $[\sigma,\pi]:=\{\lambda\in P\,|\,\sigma\le\lambda\le\pi\}$. We denote the interior of an interval by $(\sigma,\pi)$ and the half open intervals by $[\sigma,\pi)=[\sigma,\pi]\setminus\{\pi\}$ and $(\sigma,\pi]=[\sigma,\pi]\setminus\{\sigma\}$.
\end{defn}

We also recall some general poset terminology that is used throughout.

\begin{defn}
Let $P$ be a bounded poset. A \emph{chain} $c$, of length $|c|=\ell-1$, in~$P$ is a totally ordered subset of elements $c_1<c_2<\cdots<c_{\ell}$. The \emph{rank}  of an element~$\alpha\in P$, denoted $\rk_P(\alpha)$ or  $\rk(\alpha)$ when the poset is clear, is the length of the longest chain from $\hat{0}$ to $\alpha$. The rank of $P$ is given by $\rk(P)=\rk(\hat{1})$. A poset is \emph{pure} if all the maximal chains have the same length.

The \emph{join}, if it exists, of any two elements $\alpha$ and $\beta$ of a poset is the smallest element that lies above both $\alpha$ and $\beta$, and is denoted $\alpha\vee\beta$. An element $\alpha$ is \emph{covered} by $\beta$, which we denote by $\alpha\lessdot\beta$, if $\alpha<\beta$ and there is no $\kappa$ such that $\alpha<\kappa<\beta$.
\end{defn}

In this paper all posets are assumed to be pure. Two questions often asked of any poset are ``What is the M\"obius function?" and ``Is it shellable?". A poset is shellable if the maximal chains can be ordered in a ``nice" way; see \cite{Wac07} for a formal definition. Shellability has many interesting consequences for the structure and topology of a poset. The M\"obius function is defined as follows:

\begin{defn}
The M\"obius function on a poset $P$ is defined recursively, where for any elements ${a,b\in P}$ we have $\mu(a,a)=1$, $\mu(a,b)=0$ if $a\not\le b$ and if $a<b$ then:
$$\mu(a,b)=-\sum_{c\in[a,b)}\mu(a,c).$$
The \emph{M\"obius function} of bounded poset $P$ is $\mu(P)=\mu(\hat{0},\hat{1})$ and the \emph{M\"obius number} of a poset $P$ is~$\hat{\mu}(P):=\mu(\hat{P})$.
\end{defn}

We are also interested in looking at the structure of the poset. For example, it is known that a poset is not shellable if it has any disconnected subintervals of rank greater than $2$, where disconnected is defined by:
\begin{defn}
A bounded poset is \emph{disconnected} if the interior can be split into two disjoint sets, which we call \emph{components}, such that every pair of elements from separate components are incomparable.
\end{defn}

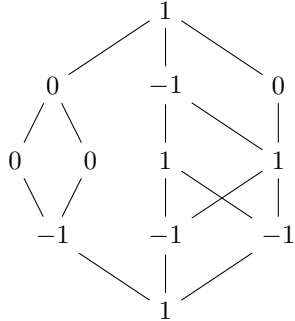
\begin{figure}\centering
\begin{tikzpicture}
\node (0) at (1.5,0) {$1$};
\node (a1) at (0,1) {$-1$};
\node (a2) at (1.5,1) {$-1$};
\node (a3) at (3,1) {$-1$};
\node (b1) at (-.5,2) {$0$};
\node (b2) at (.5,2) {$0$};
\node (b3) at (1.5,2) {$1$};
\node (b4) at (3,2) {$1$};
\node (c1) at (0,3) {$0$};
\node (c2) at (1.5,3) {$-1$};
\node (c3) at (3,3) {$0$};
\node (top) at (1.5,4) {$1$};
\draw (0) -- (a1);
\draw (a1) -- (b1);
\draw (a1) -- (b2);
\draw (b1) -- (c1);
\draw (b2) -- (c1);
\draw (c1) -- (top);
\draw (0) -- (a2);
\draw (0) -- (a3);
\draw (a2) -- (b3);
\draw (a2) -- (b4);
\draw (a3) -- (b3);
\draw (a3) -- (b4);
\draw (b3) -- (c2);
\draw (b4) -- (c2);
\draw (b4) -- (c3);
\draw (c2) -- (top);
\draw (c3) -- (top);
\end{tikzpicture}
\caption{The Hasse diagram of a disconnected bounded pure poset $P$ with rank $\rk(P)=4$ and M\"obius function $\mu(P)=1$, where each element $\alpha$ is labelled by the value of $\mu(\hat{0},\alpha)$.}\label{fig:poset}
\end{figure}

See \cref{fig:poset} for an example of a disconnected bounded pure poset. In order to study these properties of a poset we use poset fibrations, which where first introduced by Quillen in \cite{Qui78} and have many nice properties; see \cite{Bjo05} for a good overview.

\begin{defn}
A \emph{poset fibration} is a rank and order preserving surjective map between posets.
\end{defn}

Much of the terminology used here comes from poset topology, where we represent a poset $P$ as a simplicial complex $\Delta(P)$ whose faces are the chains of~$P$. Many of the properties we have introduced translate to
 topological properties on this simplicial complex. For example, $P$ is disconnected if and only if $\Delta(P)$ is disconnected, a poset fibration $f:P\rightarrow Q$ gives a fibration between~$\Delta(P)$ and~$\Delta(Q)$,
 and the M\"obius function of $P$ equals the Euler Characteristic of~$\Delta(P)$. We refer the reader to \cite{Wac07} for further background on poset topology.

\section{Pattern Posets}\label{sec:patPos}
Pattern occurrence has been studied on a variety of different combinatorial objects, with permutation patterns receiving the most attention. In this section we present pattern posets, which is a general definition of posets of combinatorial objects where the partial order is given by the occurrence of patterns. First we need to formally define pattern occurrence. 

\begin{defn}\label{defn:patRel}
We define a \emph{pattern relation}~$\sim_\rho$ as a binary relation between words on an alphabet $\Sigma$.
\end{defn}

\begin{defn}\label{defn:posCon}
Consider a word $\beta$ on an alphabet $\Sigma$ and a set of positive integers $I$. Define a \emph{position condition function} $\kappa:(\beta,I)\mapsto K$ where $K$ is a set of inequalities on $\{\beta_j|j\in I\}\cup I$.
\end{defn}

\begin{defn}\label{defn:patPair}
Consider a pair of words $\alpha,\beta$ on an alphabet $\Sigma$. Given a pattern relation $\sim_\rho$ and a position condition function $\kappa$, we call $(\rho,\kappa)$ a \emph{pattern pair}. We say that a subsequence $\beta_{i_1}\beta_{i_2}\dots\beta_{i_k}$ is a \emph{($\rho,\kappa)$-occurrence} of $\alpha$ if $\alpha\sim_\rho\beta_{i_1}\ldots\beta_{i_k}$ and $\{\beta_j|j\in I\}\cup I$ satisfies the conditions of $\kappa(\beta,I)$, where $I=\{i_1,i_2,\ldots,i_k\}$.
\end{defn}

We can use our notion of pattern relations to define a poset as follows:

\begin{defn}\label{defn:patPos}
Consider a pattern pair $(\rho,\kappa)$. Given $\sigma,\pi\in\Sigma^*$ we define a binary relation by $\sigma\le_{\rho,\kappa}\pi$ if there is a $(\rho,\kappa)$-occurrence of $\sigma$ in~$\pi$. If $\le_{\rho,\kappa}$ is reflexive, antisymmetric and transitive, then we define a \emph{pattern poset} $P(B,\rho,\kappa)$ as the poset with elements $B\subseteq\Sigma^*$ and partial order $\le_{\rho,\kappa}$.
\end{defn}

\begin{ex}\label{ex:patRel}
Let $\mathcal{S}$ and $\mathcal{A}$ be the set of all permutations and set of all words on the non-negative integers, respectively, and let $\kappa_\emptyset$ be the position condition function that always maps to the empty set. Some examples of pattern posets are:
\begin{enumerate}
\item On $\mathcal{A}$ define $\alpha\sim_\omega\beta$ if and only if $\alpha$ and $\beta$ are equal. Then $P(\mathcal{A},\omega,\kappa_\emptyset)$ is the poset of subword order; see \cite{Bjo90},
\item On $\mathcal{A}$, given a poset $Q$ whose elements are $\mathbb{N}$, define $\sigma\sim_{\omega(Q)}\pi$ if and only if $\sigma_i\le_Q\pi_i$ for all $i$. Then $P(\mathcal{A},\omega(Q),\kappa_\emptyset)$ is the poset of generalised subword order; see \cite{McnSag12},
\item Let $\mathcal{N}$ be the chain of natural numbers. Then $P(\mathcal{A},\omega(\mathcal{N}),\kappa_\emptyset)$ is the composition poset; see \cite{SagVat06},
\item On $\mathcal{S}$ define $\sigma\sim_\delta\pi$ if and only if $\red(\sigma)=\red(\pi)$. Then $P(\mathcal{S},\delta,\kappa_\emptyset)$ is the classical permutation pattern poset; see \cite{McSt13},
\item Define $\kappa_c(\beta,\mathbb{N})=\{\beta_i+1=\beta_{i+1}\,|\,\forall i\in\mathbb{N}\}$. Then $P(\mathcal{S},\delta,\kappa_c)$ is the consecutive permutation pattern poset; see \cite{BFS11,EM15},
\item Define $\kappa_X(\beta,\mathbb{N})=\{\beta_i+1=\beta_{i+1}\,|\,X_i=1\}$, where $X$ is an infinite binary string. Then $P(\mathcal{A},\omega,\kappa_X)$ are the vincular pattern posets; see \cite{BF14},
\item Define $B\subseteq\mathcal{A}$ as the set of Dyck words. Then $P(B,\omega,\kappa_\emptyset)$ is the Dyck pattern poset; see \cite{BLPW13,Bac14}.
\end{enumerate}
\end{ex}

\begin{figure}
\begin{center}\begin{tikzpicture}[scale=1]\def\x{0}\def\y{-5}\def\b{.075}
\node[label={[font=\scriptsize]left:$21$}] (21) at (0,-1) {};
\node[label={[font=\scriptsize]left:$312$}] (312) at (-1.5,0) {};
\node[label={[font=\scriptsize]left:$213$}] (213) at (-.3,0) {};
\node[label={[font=\scriptsize]right:$231$}] (231) at (.3,0) {};
\node[label={[font=\scriptsize]right:$132$}] (132) at (1.5,0) {};
\node[label={[font=\scriptsize]left:$4123$}] (4123) at (-1.5,1) {};
\node[label={[font=\scriptsize]left:$3124$}] (3124) at (-.3,1) {};
\node[label={[font=\scriptsize]right:$3142$}] (3142) at (.3,1) {};
\node[label={[font=\scriptsize]right:$1243$}] (1243) at (1.5,1) {};
\node[label={[font=\scriptsize]left:$41253$}] (41253) at (0,2) {};
\node (lab) at (0,-1.5){Classical Permutation Poset};
\draw (21) -- (312) -- (4123) -- (41253) -- (3124) -- (312) -- (3142) -- (41253) -- (1243) -- (132) -- (21) -- (231) -- (3142) -- (213) -- (21);
\draw (3124) -- (213);\draw (132) -- (3142);
\shade[shading=ball] (21) circle (\b);
\shade[shading=ball] (312) circle (\b);
\shade[shading=ball] (132) circle (\b);
\shade[shading=ball] (213) circle (\b);
\shade[shading=ball] (231) circle (\b);
\shade[shading=ball] (3124) circle (\b);
\shade[shading=ball] (1243) circle (\b);
\shade[shading=ball] (4123) circle (\b);
\shade[shading=ball] (3142) circle (\b);
\shade[shading=ball] (41253) circle (\b);
\def\x{6}
\node[label={[font=\scriptsize]left:$21$}] (21) at (0+\x,-1) {};
\node[label={[font=\scriptsize]left:$312$}] (312) at (-.5+\x,0) {};
\node[label={[font=\scriptsize]right:$132$}] (132) at (.5+\x,0) {};
\node[label={[font=\scriptsize]left:$3124$}] (3124) at (-.5+\x,1) {};
\node[label={[font=\scriptsize]right:$1243$}] (1243) at (.5+\x,1) {};
\node[label={[font=\scriptsize]left:$41253$}] (41253) at (0+\x,2) {};
\node (lab) at (0+\x,-1.5){Consecutive Permutation Poset};
\draw (21) -- (312) -- (3124) --(41253) -- (1243) -- (132) -- (21);
\shade[shading=ball] (21) circle (\b);
\shade[shading=ball] (312) circle (\b);
\shade[shading=ball] (132) circle (\b);
\shade[shading=ball] (3124) circle (\b);
\shade[shading=ball] (1243) circle (\b);
\shade[shading=ball] (41253) circle (\b);
\node[label={[font=\scriptsize]left:$2122$}] (2122) at (0,-1+\y) {};
\node[label={[font=\scriptsize]right:$21232$}] (21232) at (1.5,0+\y) {};
\node[label={[font=\scriptsize]left:$21322$}] (21322) at (0,0+\y) {};
\node[label={[font=\scriptsize]left:$23122$}] (23122) at (-1.5,0+\y) {};
\node[label={[font=\scriptsize]left:$231322$}] (231322) at (-1.5,1+\y) {};
\node[label={[font=\scriptsize]left:$231232$}] (231232) at (0,1+\y) {};
\node[label={[font=\scriptsize]right:$213232$}] (213232) at (1.5,1+\y) {};
\node[label={[font=\scriptsize]left:$2313232$}] (2313232) at (0,2+\y) {};
\node (lab) at (0,-1.5+\y){Subword Poset};
\draw (21322) -- (2122) -- (21232);\draw (2122) -- (23122);
\draw (231232) -- (2313232) -- (213232);\draw (2313232) -- (231322);
\draw (21232) -- (231232) -- (23122) -- (231322) -- (21322) -- (213232) -- (21232);
\shade[shading=ball] (2122) circle (\b);
\shade[shading=ball] (21232) circle (\b);
\shade[shading=ball] (21322) circle (\b);
\shade[shading=ball] (23122) circle (\b);
\shade[shading=ball] (231322) circle (\b);
\shade[shading=ball] (231232) circle (\b);
\shade[shading=ball] (213232) circle (\b);
\shade[shading=ball] (2313232) circle (\b);
\node[label={[font=\scriptsize]left:$21$}] (21) at (0+\x,-1+\y) {};;
\node[label={[font=\scriptsize]right:$31$}] (31) at (.5+\x,0+\y) {};
\node[label={[font=\scriptsize]left:$22$}] (22) at (-.5+\x,0+\y) {};
\node[label={[font=\scriptsize]left:$23$}] (23) at (-.5+\x,1+\y) {};
\node[label={[font=\scriptsize]right:$32$}] (32) at (.5+\x,1+\y) {};
\node[label={[font=\scriptsize]left:$33$}] (33) at (0+\x,2+\y) {};
\node (lab) at (0+\x,-1.5+\y){Composition Poset};
\draw (21) -- (31) -- (32) -- (33) -- (23) -- (22) -- (21);
\draw (22) -- (32);
\shade[shading=ball] (21) circle (\b);
\shade[shading=ball] (31) circle (\b);
\shade[shading=ball] (22) circle (\b);
\shade[shading=ball] (23) circle (\b);
\shade[shading=ball] (32) circle (\b);
\shade[shading=ball] (33) circle (\b);
\node[label={[font=\scriptsize]left:$01$}] (21) at (0+.5*\x,-1+2*\y) {};;
\node[label={[font=\scriptsize]right:$0011$}] (31) at (.5+\x*0.5,0+2*\y) {};
\node[label={[font=\scriptsize]left:$0101$}] (22) at (-.5+\x*0.5,0+2*\y) {};
\node[label={[font=\scriptsize]left:$010101$}] (23) at (-.5+\x*0.5,1+2*\y) {};
\node[label={[font=\scriptsize]right:$001101$}] (32) at (.5+\x*0.5,1+2*\y) {};
\node[label={[font=\scriptsize]left:$00110101$}] (33) at (0+\x*0.5,2+2*\y) {};
\node (lab) at (0+\x*0.5,-1.5+2*\y){Dyck Path Poset};
\draw (21) -- (31) -- (32) -- (33) -- (23) -- (22) -- (21);
\draw (22) -- (32);
\shade[shading=ball] (21) circle (\b);
\shade[shading=ball] (31) circle (\b);
\shade[shading=ball] (22) circle (\b);
\shade[shading=ball] (23) circle (\b);
\shade[shading=ball] (32) circle (\b);
\shade[shading=ball] (33) circle (\b);
\end{tikzpicture}\end{center}
\caption{Intervals of five different pattern posets.}\label{fig:PP}
\end{figure}

See \cref{fig:PP} for examples of some of the pattern posets listed above. When it is clear what poset we are considering we drop the subscript and use the notation $\le$. The notion of pattern pair is very general and not every pattern pair will induce a partial order, in fact most pattern pairs do not.

\section{A Poset Fibration}\label{sec:posfib}
In this section we introduce a poset fibration for pattern posets, along with a variation for pattern posets with a particular property that we define. Poset fibrations were first studied by Quillen in \cite{Qui78} and have many nice properties. First we introduce embeddings which play an important role throughout the paper:

\begin{defn}
Consider a pattern poset $P$ and two elements $\alpha\le\beta$ of $P$. An \emph{embedding} $\eta$ of~$\alpha$ in~$\beta$ is a sequence of length $|\beta|$, consisting of the letters of $\alpha$ and dashes, such that the non-dash letters are exactly the positions of an occurrence of $\alpha$ in $\beta$ and the removal of the dashed letters results in a word which is equal to $\alpha$. Define~$E^{\alpha,\beta}$ as the set of embeddings of $\alpha$ in $\beta$.
\end{defn}
\begin{ex}
In the composition poset $121--$ is an embedding of $121$ in~$13211$. In the classical permutation poset $-56--3$ is an embedding of $231$ in $156243$, because $563$ and $231$ are equal as they have the same reduced form. 
\end{ex}

\begin{defn}
Given an embedding $\eta\in E^{\alpha,\beta}$, we call the positions of the dashed letters in $\eta$ the \emph{empty positions} and let the \emph{zero set} of $\eta$, denoted $Z(\eta)$, be the set of empty positions in~$\eta$.
\end{defn}

Traditionally zeroes are used in embeddings instead of dashes, but we use dashes as this allows us to consider words that contain zeroes. We can define a poset fibration using the embeddings of a pattern poset:

\begin{defn}\label{defn:PF}
Consider an interval $[\sigma,\pi]$ of a pattern poset $P$. Define the poset $$A(\sigma,\pi)=\bigcup_{\lambda\in(\sigma,\pi)}E^{\lambda,\pi},$$ with the partial order $\eta\le\phi$ if~$Z(\eta)\supseteq Z(\phi)$ and~$\alpha\le_P \beta$, where $\alpha$ and $\beta$ are obtained by removing the dashes from $\eta$ and $\phi$, respectively. Also, define~$A^*(\sigma,\pi)=A(\sigma,\pi)\cup E^{\sigma,\pi}$. Moreover, define $\hat{A}(\sigma,\pi)$ and $\hat{A}^*(\sigma,\pi)$ as the bounded posets obtained by adding a top element $\pi$ and bottom element $\hat{0}$ to~$A(\sigma,\pi)$ and $A^*(\sigma,\pi)$.
\end{defn}
\begin{ex}
Consider the permutation pattern poset with $\pi=243516$, then $-435--\le-435-6$ however $-4--16\not\le-435-6$.
\end{ex}

We can now define the poset fibration:
\begin{defn}
 Define~${f^P_\pi:A(\sigma,\pi)\rightarrow (\sigma,\pi)}$ as the map that takes the elements of $E^{\lambda,\pi}$ to $\lambda$.
\end{defn}

So we have a poset fibration where $f^P_\pi$ is the projection map and $A(\sigma,\pi)$ is the total space. Where it is clear we use the notation $f$, dropping the $P$ and $\pi$. We can define a variation of this poset fibration by considering pattern posets with the following property:
\begin{defn}
Consider a pattern poset $P$ whose elements are defined on the alphabet $\Sigma$. We say that $P$ is \emph{closed} if we can represent $\Sigma\cup\{-\}$ as a tree rooted at $-$, which we denote $\bar{\Sigma}$, where for every $\eta\in \hat{A}^*(\sigma,\pi)\setminus\{\hat{0}\}$ we have:
\begin{equation}\label{eq:closed}[\eta,\pi]=\{\phi\,|\,\eta_i\le_{\bar{\Sigma}}\phi_i\le_{\bar{\Sigma}}\pi_i\text{ for all }i\}\end{equation}
for every pair $\sigma,\pi\in P$. Moreover, we say $P$ is \emph{fully-closed} if $\bar{\Sigma}$ consists of an antichain and the bottom element $-$.
\end{defn}
\begin{defn}
In a closed pattern poset $P$ we say that we are \emph{decreasing} a letter in an embedding if we change a letter $i$ to $j\lessdot_{\bar{\Sigma}} i$, and \emph{increasing} if we change $j$ to $i$. Similarly, we can decrease a letter in an element $\pi\in P$, and if we reach $-$ then we delete the letter.
\end{defn}
\begin{ex}
The subword order and classical permutation pattern posets are fully-closed, because we can construct $\bar{\Sigma}$ from the antichain $\mathbb{N}$ and bottom element $-$, and we can get between embeddings by turning non-empty letter to empty letters, and vice-versa.

The composition poset is closed but not fully-closed, where $\bar{\Sigma}$ is the chain of natural numbers, and $-$ at the bottom.

The consecutive permutation poset is not closed, because only embeddings where the non-empty letters are consecutive are valid embeddings, so we can never satisfy \eqref{eq:closed}. Similarly, the Dyck pattern poset is not closed, because embeddings must have an even number of non-empty letters.
\end{ex}

\begin{rem}\label{rem:closed}
If $[\sigma,\pi]$ is an interval of a closed pattern poset $P$, then every interval $[\eta,\phi]$ in $\hat{A}^*(\sigma,\pi)$, with $\eta\not=\hat{0}$, is isomorphic to the poset $\prod_i[\eta_i,\phi_i]$, where $[\eta_i,\phi_i]$ is a chain in $\bar{\Sigma}$. Moreover, if $P$ is fully closed then $[\eta,\phi]$ is isomorphic to the boolean lattice $B_t$, where $t$ is the number of positions $\eta_i\not=\phi_i$. This is because a chain $[\eta_i,\phi_i]$ is of length $1$ if $\eta_i\not=\phi_i$, and of length $0$ otherwise. 
\end{rem}

As $\bar{\Sigma}$ is a tree the decreasing and increasing operations are well defined on any embedding $\eta\in E^{\sigma,\pi}$ as there is a unique path between $-$ and $\pi_i$ in $\bar{\Sigma}$, for all $i$. The decreasing operation on a word $\pi$ is well defined as this corresponds to travelling down the unique path from $\pi_i$ to $-$. However, we cannot always uniquely increase a letter of a word, because an element of $\bar{\Sigma}$ can be covered by more than one element.

Next we introduce normal embeddings, which haved played an important role in many of the existing results on the M\"obius function of pattern posets. For example, normal embeddings appear in results  on the poset of words with subword order \cite{Bjo90,Bjo93}, the poset of words with composition order \cite{SagVat06}, the poset of words with generalised subword order~\cite{SagVat06,McnSag12} and the classical permutation poset \cite{BJJS11,Smith13,Smith14,Smith15}. We generalise these notions of normal embeddings and present a simple definition for a normal embedding for any closed pattern poset. This definition is different from some definitions in the literature but extracts the common aspect of all of them, and the variation is then accounted for in our formula for the M\"obius function in \cref{sec:results}. First we introduce \emph{adjacencies} which play an important role in defining a normal embedding.

\begin{defn}
Consider a closed pattern poset $P:=P(B,\rho,\kappa)$. An \emph{adjacency} in an element $\sigma\in P$ is a maximal sequence of consecutive positions such that decreasing any letter of the adjacency yields the same element relative to $\sim_{\rho}$. An adjacency of length $1$ is \emph{trivial}, the \emph{tail} of a non-trivial adjacency is all but the first letter of the adjacency and trivial adjacencies have no tails.
\end{defn}
\begin{ex}
Consider the classical permutation poset, where the decreasing operation is deletion. In the permutation $\pi=2341657$ we have an adjacency $234$ because deletion of any letter gives the permutation $231567$. So the adjacencies are $234,\,1,\,65,\,7$ and the tails are $34$ and $5$.

Consider the composition poset, where the decreasing operation is: reduce the value of $i\ge2$ by $1$ or delete $1$. In $321122$ the only non-trivial adjacency is $11$. In fact, in any word of this poset the only non-trivial adjacencies are consecutive sequences of $1$'s.
\end{ex}

\begin{defn}
Given any embedding $\eta$ of $\sigma$ in $\pi$ in a closed pattern poset we say a position $i$ of $\eta$ is \emph{full} if $\eta_i=\pi_i$ and \emph{fillable} if increasing $\eta_i$ once results in $\pi_i$.
\end{defn}

For example, a position is full in a fully-closed pattern poset if and only if it is non-empty and fillable if and only if it is empty. A position is full in the composition poset if $\eta_i=\pi_i$ and fillable if~$\eta_i=\pi_i-1$ or $\eta_i=-$ and~$\pi_i=1$. 

Using our definition of adjacency we can define a normal embedding, which appears frequently in the results on pattern posets:

\begin{defn}\label{defn:normal}
Consider a closed pattern poset $P:=P(B,\rho,\kappa)$ and two elements $\sigma,\pi\in P$. An embedding~$\eta$ of $\sigma$ in $\pi$  is \emph{normal} if all the positions that are in a tail of any adjacency in $\pi$ are full in $\eta$ and all other positions are fillable. Let $\NE(\sigma,\pi)$ denote the number of normal embeddings of $\sigma$ in $\pi$.

An embedding $\eta$ is \emph{representative} if for every adjacency of $\pi$ the corresponding letters in $\eta$ have all the empty letters to the left, the full letters to the right and at most one non-full non-empty letter positioned between them. Let $\hat{E}^{\sigma,\pi}$ be the set of representative embeddings of $\sigma$ in $\pi$.
\end{defn}
\begin{ex}
Consider the classical permutation poset. The embeddings of~$213$ in $231645$ are:
$$2-16--\,\,\,\,\,\,\,\,-316--\,\,\,\,\,\,\,\,2-1-4-\,\,\,\,\,\,\,\,-31-4-\,\,\,\,\,\,\,\,2-1--5\,\,\,\,\,\,\,\,-31--5$$
The representative embeddings are $-316--$ and $-31--5$ and the only normal embedding is~$-31--5$.
\end{ex}
Note that in a fully-closed pattern poset an embedding is representative if there is no empty position to the right of a non-empty position in the same adjacency, and normal if all positions in the tail of an adjacency are non-empty.

The definition of normal in \cref{defn:normal} is equivalent to the definition of normal in the poset of words with subword order; see \cite{Bjo90}, and the classical permutation poset; see \cite{Smith15}. However, the definition is not equivalent to the definitions given for the poset of words with composition order~\cite{SagVat06} or generalised subword order \cite{McnSag12}. In these cases our definition of normal embeddings gives a subset of the normal embeddings according to the previous definitions. We account for these differences in the formulas we present in \cref{sec:mobGen}.

Using the notion of representative embeddings we can define another poset fibration, with a smaller total space than the initial fibration. This allows us to simplify many of the results presented in \cref{sec:results}.

\begin{defn}
Consider an interval $[\sigma,\pi]$ of a closed pattern poset. Let $R(\sigma,\pi)$ and $R^*(\sigma,\pi)$ be the posets of all the representative embeddings of~$\lambda$ in~$\pi$, for all $\lambda\in(\sigma,\pi)$ and $\lambda\in[\sigma,\pi)$, respectively. Moreover, let $\hat{R}(\sigma,\pi)$ and~$\hat{R}^*(\sigma,\pi)$ be the bounded posets obtained by adding the top element $\pi$ and bottom element~$\hat{0}$.
\end{defn}

\begin{rem}\label{rem:prodChains}
Consider an interval $[\sigma,\pi]$ of a closed pattern poset and any $\eta\in R^*(\sigma,\pi)$. To obtain an element that covers $\eta$ in $\hat{R}^*(\sigma,\pi)$ we must increase the rightmost non-full position of an adjacency, as increasing any other position would result in a non-representative embedding. So let $\hat{\eta}_i$ be the number of increasing operations required until every letter of the $i$'th adjacency is full and let~$[0,\hat{\eta}_i]$ be the chain of integers from $0$ to $\hat{\eta}_i$. The interval $[\eta,\pi]$ in $\hat{R}^*(\sigma,\pi)$ is isomorphic to the product of chains~$[0,\hat{\eta}_1]\times\cdots\times[0,\hat{\eta}_t]$.
\end{rem}

The poset of representative embeddings $R(\sigma,\pi)$ is a subposet of the poset of all embedding $A(\sigma,\pi)$. So we can define a poset fibration onto $[\sigma,\pi]$ by restricting~$f^P_\pi$ to $R(\sigma,\pi)$, which has the total space $R(\sigma,\pi)$ and the projection map~$f^P_\pi|_{R(\sigma,\pi)}$. When the context is clear we simply use $f$ to denote the projection map.

\section{Results on Pattern Poset}\label{sec:results}

\subsection{The M\"obius Function of Intervals of a Pattern Poset}\label{sec:mobGen}
In this subsection we focus on the M\"obius function of pattern posets. The following result, which is the dual of Corollary 3.2 in~\cite{Wal81}, proves very useful:
\begin{prop}\label{prop:wal81}
Given a poset fibration $f:P\rightarrow Q$:
$$\hat{\mu}(Q)=\hat{\mu}(P)+\sum_{q\in Q}\hat{\mu}(Q_{<q})\hat{\mu}(f^{-1}(Q_{\ge q})).$$
\end{prop}

In many of our applications of \cref{prop:wal81} we consider bounded posets so use $\mu$ rather than $\hat{\mu}$. Applying \cref{prop:wal81} to the poset fibrations given in \cref{sec:posfib} gives the following results:

\begin{thm}\label{cor:mobGen}
If  $[\sigma,\pi]$ is an interval of a pattern poset, then:
\begin{align}\mu(\sigma,\pi)&=\mu(\hat{A}(\sigma,\pi))+\sum_{\lambda\in(\sigma,\pi)}\mu(\sigma,\lambda)\mu(\hat{A}^*(\lambda,\pi))\label{eq:muA0}\\&=\sum_{\eta\in E^{\sigma,\pi}}\mu(\eta,\pi)+\sum_{\lambda\in[\sigma,\pi)}\mu(\sigma,\lambda)\mu(\hat{A}^*(\lambda,\pi))\label{eq:muA4}.\end{align}
\begin{proof}
Applying \cref{prop:wal81} to the poset fibration given in \cref{defn:PF} gives \cref{eq:muA0}. The posets $A(\sigma,\pi)$ and $A^*(\sigma,\pi)$ can be considered as the union of the intervals~$(\eta,\pi)$ and $[\eta,\pi)$, respectively, for all $\eta\in E^{\sigma,\pi}$. Applying an inclusion-exclusion argument for the M\"obius function we~get:
\begin{equation}\label{eq:muA}\mu(\hat{A}(\sigma,\pi))=\sum_{\eta\in E^{\sigma,\pi}}\mu(\eta,\pi)+\sum_{\substack{S\subseteq E^{\sigma,\pi}\\|S|>1}}(-1)^{|S|}\hat{\mu}\left(\bigcap_{\eta\in S}(\eta,\pi)\right),\end{equation}
\begin{equation}\label{eq:muA3}\mu(\hat{A}^*(\sigma,\pi))=\sum_{\eta\in E^{\sigma,\pi}}\hat{\mu}([\eta,\pi))+\sum_{\substack{S\subseteq E^{\sigma,\pi}\\|S|>1}}(-1)^{|S|}\hat{\mu}\left(\bigcap_{\eta\in S}(\eta,\pi)\right).\end{equation}

Note that in \cref{eq:muA3} we use the intersections of $(\eta,\pi)$ instead of~$[\eta,\pi)$, because these are equivalent as $\eta$ will never be in the intersections. Moreover,~$[\eta,\pi)$ has the unique bottom element~$\eta$ and thus the M\"obius number equals $0$. So the first term on the right hand side of \cref{eq:muA3} equals zero. Therefore, the second term on the right hand side of \cref{eq:muA} is equal to $\mu(\hat{A}^*(\sigma,\pi))$, so we have:
\begin{equation}\label{eq:muA2}\mu(\hat{A}(\sigma,\pi))=\sum_{\eta\in E^{\sigma,\pi}}\mu(\eta,\pi)+\mu(\hat{A}^*(\sigma,\pi)).\end{equation}
Combining Equations \eqref{eq:muA0} and \eqref{eq:muA2} gives \cref{eq:muA4}.
\end{proof}
\end{thm}

\begin{thm}\label{thm:muEta1}
If $[\sigma,\pi]$ is an interval of a fully-closed pattern poset, then:
\begin{equation}\label{eq:main1}\mu(\sigma,\pi)=(-1)^{|\pi|-|\sigma|}E(\sigma,\pi)+\sum_{\lambda\in[\sigma,\pi)}\mu(\sigma,\lambda)\mu(\hat{A}^*(\lambda,\pi)).\end{equation}
\begin{proof}

By \cref{rem:closed} we know that in a fully-closed pattern poset the interval~$[\eta,\pi]$ is a boolean lattice of rank $|\pi|-|\sigma|$, for all $\eta\in E^{\sigma,\pi}$, so has M\"obius number~$(-1)^{|\pi|-|\sigma|}$. Therefore, the result follows from \cref{cor:mobGen}.
\end{proof}
\end{thm}

\begin{thm}\label{thm:muEta2}
If $[\sigma,\pi]$ is an interval of a closed pattern poset, then:
\begin{equation}\label{eq:main2}\mu(\sigma,\pi)=(-1)^{|\pi|-|\sigma|}\NE(\sigma,\pi)+\sum_{\lambda\in[\sigma,\pi)}\mu(\sigma,\lambda)\mu(\hat{R}^*(\lambda,\pi)).\end{equation}
\begin{proof}
If we consider the posets $R(\sigma,\pi)$ and $R^*(\sigma,\pi)$ and apply an analogous argument to that used in the proof of \cref{cor:mobGen} to derive \cref{eq:muA4}, then we get the following equation:
\begin{equation}\label{eq:muA5}\mu(\sigma,\pi)=\sum_{\eta\in \hat{E}^{\sigma,\pi}}\mu(\eta,\pi)+\sum_{\lambda\in[\sigma,\pi)}\mu(\sigma,\lambda)\mu(\hat{R}^*(\lambda,\pi))\end{equation}

By \cref{rem:prodChains} we know that $[\eta,\pi]$ is the Cartesian product of chains, for any~$\eta\in\hat{E}^{\sigma,\pi}$. Moreover, these chains all have length at most $1$ if and only if~$\eta$ is normal. If $\eta$ is normal then there are $|\pi|-|\sigma|$ chains of length $1$, the rest having length~$0$. The M\"obius function of a chain is $1$ if the chain has length~$0$,~$-1$ if the chain has length~$1$ and $0$ otherwise. Furthermore, the M\"obius function of the Cartesian product of posets is the product of the M\"obius functions. Therefore,~$\mu(\eta,\pi)$ equals~$(-1)^{|\pi|-|\sigma|}$ if~$\eta$ is normal and $0$ otherwise. So the first term on the right hand side of \cref{eq:muA5} equals~$(-1)^{|\pi|-|\sigma|}\NE(\sigma,\pi)$, which completes the proof.
\end{proof}
\end{thm}

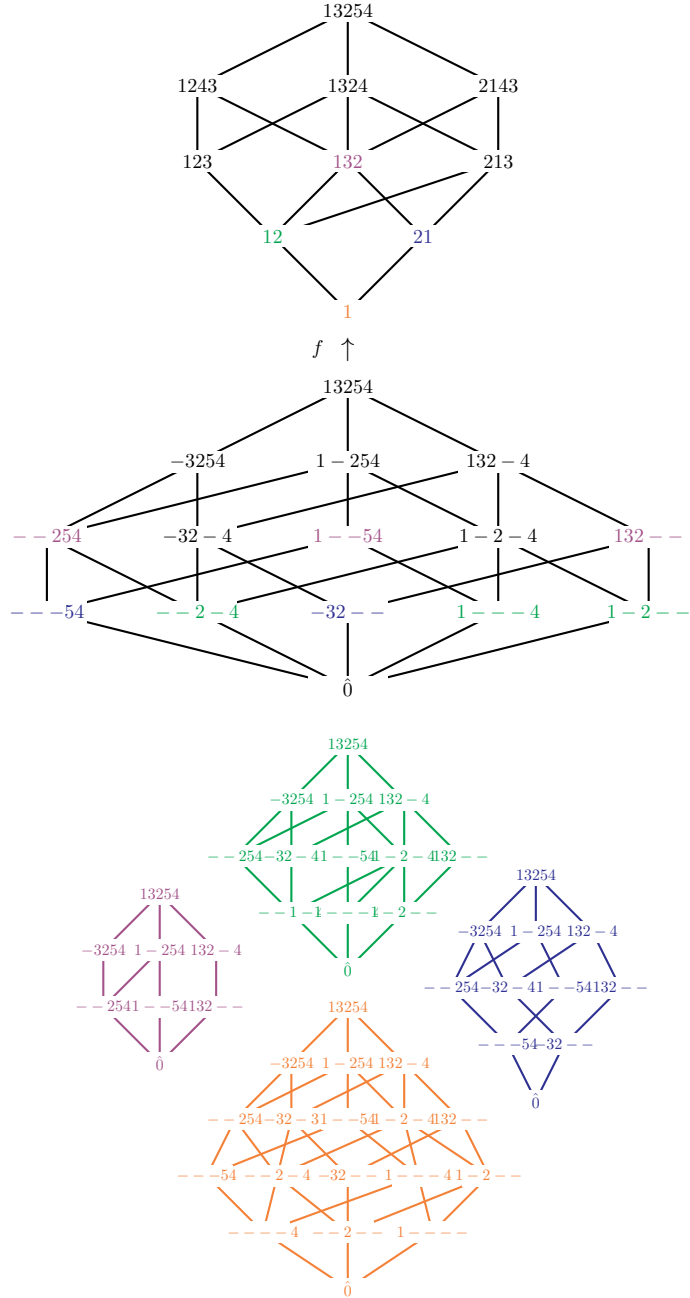
\begin{figure}\centering
\begin{tikzpicture}[scale=0.5]
\def\font{.75}
\def\bwid{30}
\def\bhei{8}
\def\x{4}
\def\y{2}
\def\d{0}
\def\dd{-10}
\node (13254) at (2*\x,4*\y+\d) {};
\node (1243) at (1*\x,3*\y+\d) {};
\node (1324) at (2*\x,3*\y+\d) {};
\node (2143) at (3*\x,3*\y+\d) {};
\node (123) at (1*\x,2*\y+\d) {};
\node (132) at (2*\x,2*\y+\d) {};
\node (213) at (3*\x,2*\y+\d) {};
\node (12) at (1.5*\x,1*\y+\d) {};
\node (21) at (2.5*\x,1*\y+\d) {};
\node (1) at (2*\x,0*\y+\d) {};
\draw[thick] (1) -- (12) -- (123) -- (1243) -- (13254);
\draw[thick] (1) -- (21) -- (132) -- (1324) -- (13254);
\draw[thick] (12) -- (132) -- (2143) -- (13254);
\draw[thick] (12) -- (213) -- (1324);
\draw[thick] (21) -- (213) -- (2143);
\draw[thick] (123) -- (1324);
\draw[thick] (132) -- (1243);
\node[label=center:\scalebox{\font}{$13254$},rectangle,inner sep=0,minimum width=\bwid pt,minimum height=\bhei pt,fill=white] at (13254){};
\node[label=center:\scalebox{\font}{$1243$},rectangle,inner sep=0,minimum width=\bwid pt,minimum height=\bhei pt,fill=white] at (1243){};
\node[label=center:\scalebox{\font}{$1324$},rectangle,inner sep=0,minimum width=\bwid pt,minimum height=\bhei pt,fill=white] at (1324){};
\node[label=center:\scalebox{\font}{$2143$},rectangle,inner sep=0,minimum width=\bwid pt,minimum height=\bhei pt,fill=white] at (2143){};
\node[label=center:\scalebox{\font}{$123$},rectangle,inner sep=0,minimum width=\bwid pt,minimum height=\bhei pt,fill=white] at (123){};
\node[label=center:\scalebox{\font}{\textcolor{DarkOrchid}{$132$}},rectangle,inner sep=0,minimum width=\bwid pt,minimum height=\bhei pt,fill=white] at (132){};
\node[label=center:\scalebox{\font}{$213$},rectangle,inner sep=0,minimum width=\bwid pt,minimum height=\bhei pt,fill=white] at (213){};
\node[label=center:\scalebox{\font}{\textcolor{Green}{$12$}},rectangle,inner sep=0,minimum width=\bwid pt,minimum height=\bhei pt,fill=white] at (12){};
\node[label=center:\scalebox{\font}{\textcolor{Blue}{$21$}},rectangle,inner sep=0,minimum width=\bwid pt,minimum height=\bhei pt,fill=white] at (21){};
\node[label=center:\scalebox{\font}{\textcolor{Orange}{$1$}},rectangle,inner sep=0,minimum width=\bwid pt,minimum height=\bhei pt,fill=white] at (1){};
\node (arrow) at (2*\x,-0.5*\y){\scalebox{1}{$\uparrow$}};
\node (f) at (1.8*\x,-0.5*\y){\scalebox{0.75}{$f$}};
\node (13254) at (2*\x,4*\y+\dd) {};
\node (02143) at (1*\x,3*\y+\dd) {};
\node (10243) at (2*\x,3*\y+\dd) {};
\node (13204) at (3*\x,3*\y+\dd) {};
\node (00132) at (0*\x,2*\y+\dd) {};
\node (02103) at (1*\x,2*\y+\dd) {};
\node (10032) at (2*\x,2*\y+\dd) {};
\node (10203) at (3*\x,2*\y+\dd) {};
\node (13200) at (4*\x,2*\y+\dd) {};
\node (00021) at (0*\x,1*\y+\dd) {};
\node (00102) at (1*\x,1*\y+\dd) {};
\node (02100) at (2*\x,1*\y+\dd) {};
\node (10002) at (3*\x,1*\y+\dd) {};
\node (10200) at (4*\x,1*\y+\dd) {};
\node (0) at (2*\x,0*\y+\dd) {};
\draw[thick] (0) -- (10002) -- (10032) -- (10243) -- (13254);
\draw[thick] (0) -- (02100) -- (02103) -- (02143) -- (13254);
\draw[thick] (0) -- (00102) -- (02103) -- (13204) -- (13254);
\draw[thick] (0) -- (00021) -- (00132) -- (02143);
\draw[thick] (0) -- (10200) -- (13200) -- (13204);
\draw[thick] (10002) -- (10203) -- (13204);
\draw[thick] (00102) -- (00132) -- (10243);
\draw[thick] (10200) -- (10203) -- (10243);
\draw[thick] (00021) -- (10032);
\draw[thick] (00102) -- (10203);
\draw[thick] (02100) -- (13200);
\node[label=center:\scalebox{\font}{$13254$},rectangle,inner sep=0,minimum width=\bwid pt,minimum height=\bhei pt,fill=white] at (13254){};
\node[label=center:\scalebox{\font}{$-3254$},rectangle,inner sep=0,minimum width=\bwid pt,minimum height=\bhei pt,fill=white] at (02143){};
\node[label=center:\scalebox{\font}{$1-254$},rectangle,inner sep=0,minimum width=\bwid pt,minimum height=\bhei pt,fill=white] at (10243){};
\node[label=center:\scalebox{\font}{$132-4$},rectangle,inner sep=0,minimum width=\bwid pt,minimum height=\bhei pt,fill=white] at (13204){};
\node[label=center:\scalebox{\font}{\textcolor{DarkOrchid}{$--254$}},rectangle,inner sep=0,minimum width=\bwid pt,minimum height=\bhei pt,fill=white] at (00132){};
\node[label=center:\scalebox{\font}{$-32-4$},rectangle,inner sep=0,minimum width=\bwid pt,minimum height=\bhei pt,fill=white] at (02103){};
\node[label=center:\scalebox{\font}{\textcolor{DarkOrchid}{$1--54$}},rectangle,inner sep=0,minimum width=\bwid pt,minimum height=\bhei pt,fill=white] at (10032){};
\node[label=center:\scalebox{\font}{$1-2-4$},rectangle,inner sep=0,minimum width=\bwid pt,minimum height=\bhei pt,fill=white] at (10203){};
\node[label=center:\scalebox{\font}{\textcolor{DarkOrchid}{$132--$}},rectangle,inner sep=0,minimum width=\bwid pt,minimum height=\bhei pt,fill=white] at (13200){};
\node[label=center:\scalebox{\font}{\textcolor{Blue}{$---54$}},rectangle,inner sep=0,minimum width=\bwid pt,minimum height=\bhei pt,fill=white] at (00021){};
\node[label=center:\scalebox{\font}{\textcolor{Green}{$--2-4$}},rectangle,inner sep=0,minimum width=\bwid pt,minimum height=\bhei pt,fill=white] at (00102){};
\node[label=center:\scalebox{\font}{\textcolor{Blue}{$-32--$}},rectangle,inner sep=0,minimum width=\bwid pt,minimum height=\bhei pt,fill=white] at (02100){};
\node[label=center:\scalebox{\font}{\textcolor{Green}{$1---4$}},rectangle,inner sep=0,minimum width=\bwid pt,minimum height=\bhei pt,fill=white] at (10002){};
\node[label=center:\scalebox{\font}{\textcolor{Green}{$1-2--$}},rectangle,inner sep=0,minimum width=\bwid pt,minimum height=\bhei pt,fill=white] at (10200){};
\node[label=center:\scalebox{\font}{$\hat{0}$},rectangle,inner sep=0,minimum width=\bwid pt,minimum height=\bhei pt,fill=white] at (0){};
\def\font{.6}
\def\l{10}
\def\ll{5}
\def\lll{0}
\def\x{1.5}
\def\y{1.5}
\def\ddd{-26}
\def\dddd{-21}
\def\ddddx{-20}
\def\ddddd{-17.5}
\def\bwid{20}
\def\bhei{8}
\node (13254) at (2*\x+\ll,4*\y+\ddddd) {};
\node (02143) at (1*\x+\ll,3*\y+\ddddd) {};
\node (10243) at (2*\x+\ll,3*\y+\ddddd) {};
\node (13204) at (3*\x+\ll,3*\y+\ddddd) {};
\node (00132) at (0*\x+\ll,2*\y+\ddddd) {};
\node (02103) at (1*\x+\ll,2*\y+\ddddd) {};
\node (10032) at (2*\x+\ll,2*\y+\ddddd) {};
\node (10203) at (3*\x+\ll,2*\y+\ddddd) {};
\node (13200) at (4*\x+\ll,2*\y+\ddddd) {};
\node (00102) at (1*\x+\ll,1*\y+\ddddd) {};
\node (10002) at (2*\x+\ll,1*\y+\ddddd) {};
\node (10200) at (3*\x+\ll,1*\y+\ddddd) {};
\node (0) at (2*\x+\ll,0*\y+\ddddd) {};
\draw[Green,thick] (0) -- (10002) -- (10032) -- (10243) -- (13254);
\draw[Green,thick] (02103) -- (02143) -- (13254);
\draw[Green,thick] (0) -- (00102) -- (02103) -- (13204) -- (13254);
\draw[Green,thick] (00132) -- (02143);
\draw[Green,thick] (0) -- (10200) -- (13200) -- (13204);
\draw[Green,thick] (10002) -- (10203) -- (13204);
\draw[Green,thick] (00102) -- (00132) -- (10243);
\draw[Green,thick] (10200) -- (10203) -- (10243);
\draw[Green,thick] (00102) -- (10203);
\node[label=center:\scalebox{\font}{\textcolor{Green}{$13254$}},rectangle,inner sep=0,minimum width=\bwid pt,minimum height=\bhei pt,fill=white] at (13254){};
\node[label=center:\scalebox{\font}{\textcolor{Green}{$-3254$}},rectangle,inner sep=0,minimum width=\bwid pt,minimum height=\bhei pt,fill=white] at (02143){};
\node[label=center:\scalebox{\font}{\textcolor{Green}{$1-254$}},rectangle,inner sep=0,minimum width=\bwid pt,minimum height=\bhei pt,fill=white] at (10243){};
\node[label=center:\scalebox{\font}{\textcolor{Green}{$132-4$}},rectangle,inner sep=0,minimum width=\bwid pt,minimum height=\bhei pt,fill=white] at (13204){};
\node[label=center:\scalebox{\font}{\textcolor{Green}{$--254$}},rectangle,inner sep=0,minimum width=\bwid pt,minimum height=\bhei pt,fill=white] at (00132){};
\node[label=center:\scalebox{\font}{\textcolor{Green}{$-32-4$}},rectangle,inner sep=0,minimum width=\bwid pt,minimum height=\bhei pt,fill=white] at (02103){};
\node[label=center:\scalebox{\font}{\textcolor{Green}{$1--54$}},rectangle,inner sep=0,minimum width=\bwid pt,minimum height=\bhei pt,fill=white] at (10032){};
\node[label=center:\scalebox{\font}{\textcolor{Green}{$1-2-4$}},rectangle,inner sep=0,minimum width=\bwid pt,minimum height=\bhei pt,fill=white] at (10203){};
\node[label=center:\scalebox{\font}{\textcolor{Green}{$132--$}},rectangle,inner sep=0,minimum width=\bwid pt,minimum height=\bhei pt,fill=white] at (13200){};
\node[label=center:\scalebox{\font}{\textcolor{Green}{$--1-4$}},rectangle,inner sep=0,minimum width=\bwid pt,minimum height=\bhei pt,fill=white] at (00102){};
\node[label=center:\scalebox{\font}{\textcolor{Green}{$1---4$}},rectangle,inner sep=0,minimum width=\bwid pt,minimum height=\bhei pt,fill=white] at (10002){};
\node[label=center:\scalebox{\font}{\textcolor{Green}{$1-2--$}},rectangle,inner sep=0,minimum width=\bwid pt,minimum height=\bhei pt,fill=white] at (10200){};
\node[label=center:\scalebox{\font}{\textcolor{Green}{$\hat{0}$}},rectangle,inner sep=0,minimum width=\bwid pt,minimum height=\bhei pt,fill=white] at (0){};
\node (13254) at (2*\x+\l,4*\y+\dddd) {};
\node (02143) at (1*\x+\l,3*\y+\dddd) {};
\node (10243) at (2*\x+\l,3*\y+\dddd) {};
\node (13204) at (3*\x+\l,3*\y+\dddd) {};
\node (00132) at (0.5*\x+\l,2*\y+\dddd) {};
\node (02103) at (1.5*\x+\l,2*\y+\dddd) {};
\node (10032) at (2.5*\x+\l,2*\y+\dddd) {};
\node (13200) at (3.5*\x+\l,2*\y+\dddd) {};
\node (00021) at (1.5*\x+\l,1*\y+\dddd) {};
\node (02100) at (2.5*\x+\l,1*\y+\dddd) {};
\node (0) at (2*\x+\l,0*\y+\dddd) {};
\draw[Blue,thick] (10032) -- (10243) -- (13254);
\draw[Blue,thick] (0) -- (02100) -- (02103) -- (02143) -- (13254);
\draw[Blue,thick] (02103) -- (13204) -- (13254);
\draw[Blue,thick] (0) -- (00021) -- (00132) -- (02143);
\draw[Blue,thick] (13200) -- (13204);
\draw[Blue,thick] (00132) -- (10243);
\draw[Blue,thick] (00021) -- (10032);
\draw[Blue,thick] (02100) -- (13200);
\node[label=center:\scalebox{\font}{\textcolor{Blue}{$13254$}},rectangle,inner sep=0,minimum width=\bwid pt,minimum height=\bhei pt,fill=white] at (13254){};
\node[label=center:\scalebox{\font}{\textcolor{Blue}{$-3254$}},rectangle,inner sep=0,minimum width=\bwid pt,minimum height=\bhei pt,fill=white] at (02143){};
\node[label=center:\scalebox{\font}{\textcolor{Blue}{$1-254$}},rectangle,inner sep=0,minimum width=\bwid pt,minimum height=\bhei pt,fill=white] at (10243){};
\node[label=center:\scalebox{\font}{\textcolor{Blue}{$132-4$}},rectangle,inner sep=0,minimum width=\bwid pt,minimum height=\bhei pt,fill=white] at (13204){};
\node[label=center:\scalebox{\font}{\textcolor{Blue}{$--254$}},rectangle,inner sep=0,minimum width=\bwid pt,minimum height=\bhei pt,fill=white] at (00132){};
\node[label=center:\scalebox{\font}{\textcolor{Blue}{$-32-4$}},rectangle,inner sep=0,minimum width=\bwid pt,minimum height=\bhei pt,fill=white] at (02103){};
\node[label=center:\scalebox{\font}{\textcolor{Blue}{$1--54$}},rectangle,inner sep=0,minimum width=\bwid pt,minimum height=\bhei pt,fill=white] at (10032){};
\node[label=center:\scalebox{\font}{\textcolor{Blue}{$132--$}},rectangle,inner sep=0,minimum width=\bwid pt,minimum height=\bhei pt,fill=white] at (13200){};
\node[label=center:\scalebox{\font}{\textcolor{Blue}{$---54$}},rectangle,inner sep=0,minimum width=\bwid pt,minimum height=\bhei pt,fill=white] at (00021){};
\node[label=center:\scalebox{\font}{\textcolor{Blue}{$-32--$}},rectangle,inner sep=0,minimum width=\bwid pt,minimum height=\bhei pt,fill=white] at (02100){};
\node[label=center:\scalebox{\font}{\textcolor{Blue}{$\hat{0}$}},rectangle,inner sep=0,minimum width=\bwid pt,minimum height=\bhei pt,fill=white] at (0){};
\node (13254) at (2*\x+\lll,3*\y+\ddddx) {};
\node (02143) at (1*\x+\lll,2*\y+\ddddx) {};
\node (10243) at (2*\x+\lll,2*\y+\ddddx) {};
\node (13204) at (3*\x+\lll,2*\y+\ddddx) {};
\node (00132) at (1*\x+\lll,1*\y+\ddddx) {};
\node (10032) at (2*\x+\lll,1*\y+\ddddx) {};
\node (13200) at (3*\x+\lll,1*\y+\ddddx) {};
\node (0) at (2*\x+\lll,0*\y+\ddddx) {};
\draw[DarkOrchid,thick] (0) -- (10032) -- (10243) -- (13254);
\draw[DarkOrchid,thick] (0) -- (00132) -- (02143) -- (13254);
\draw[DarkOrchid,thick] (0) -- (13200) -- (13204) -- (13254);
\draw[DarkOrchid,thick] (00132) -- (10243);
\node[label=center:\scalebox{\font}{\textcolor{DarkOrchid}{$13254$}},rectangle,inner sep=0,minimum width=\bwid pt,minimum height=\bhei pt,fill=white] at (13254){};
\node[label=center:\scalebox{\font}{\textcolor{DarkOrchid}{$-3254$}},rectangle,inner sep=0,minimum width=\bwid pt,minimum height=\bhei pt,fill=white] at (02143){};
\node[label=center:\scalebox{\font}{\textcolor{DarkOrchid}{$1-254$}},rectangle,inner sep=0,minimum width=\bwid pt,minimum height=\bhei pt,fill=white] at (10243){};
\node[label=center:\scalebox{\font}{\textcolor{DarkOrchid}{$132-4$}},rectangle,inner sep=0,minimum width=\bwid pt,minimum height=\bhei pt,fill=white] at (13204){};
\node[label=center:\scalebox{\font}{\textcolor{DarkOrchid}{$--254$}},rectangle,inner sep=0,minimum width=\bwid pt,minimum height=\bhei pt,fill=white] at (00132){};
\node[label=center:\scalebox{\font}{\textcolor{DarkOrchid}{$1--54$}},rectangle,inner sep=0,minimum width=\bwid pt,minimum height=\bhei pt,fill=white] at (10032){};
\node[label=center:\scalebox{\font}{\textcolor{DarkOrchid}{$132--$}},rectangle,inner sep=0,minimum width=\bwid pt,minimum height=\bhei pt,fill=white] at (13200){};
\node[label=center:\scalebox{\font}{\textcolor{DarkOrchid}{$\hat{0}$}},rectangle,inner sep=0,minimum width=\bwid pt,minimum height=\bhei pt,fill=white] at (0){};
\node (13254) at (2*\x+\ll,5*\y+\ddd) {};
\node (02143) at (1*\x+\ll,4*\y+\ddd) {};
\node (10243) at (2*\x+\ll,4*\y+\ddd) {};
\node (13204) at (3*\x+\ll,4*\y+\ddd) {};
\node (00132) at (0*\x+\ll,3*\y+\ddd) {};
\node (02103) at (1*\x+\ll,3*\y+\ddd) {};
\node (10032) at (2*\x+\ll,3*\y+\ddd) {};
\node (10203) at (3*\x+\ll,3*\y+\ddd) {};
\node (13200) at (4*\x+\ll,3*\y+\ddd) {};
\node (00021) at (-0.5*\x+\ll,2*\y+\ddd) {};
\node (00102) at (0.75*\x+\ll,2*\y+\ddd) {};
\node (02100) at (2*\x+\ll,2*\y+\ddd) {};
\node (10002) at (3.25*\x+\ll,2*\y+\ddd) {};
\node (10200) at (4.5*\x+\ll,2*\y+\ddd) {};
\node (00001) at (0.5*\x+\ll,1*\y+\ddd) {};
\node (00100) at (2*\x+\ll,1*\y+\ddd) {};
\node (10000) at (3.5*\x+\ll,1*\y+\ddd) {};
\node (0) at (2*\x+\ll,0*\y+\ddd) {};
\draw[Orange,thick] (0) -- (10000) -- (10002) -- (10032) -- (10243) -- (13254);
\draw[Orange,thick] (0) -- (00100) -- (02100) -- (02103) -- (02143) -- (13254);
\draw[Orange,thick] (0) -- (00001) -- (00102) -- (02103) -- (13204) -- (13254);
\draw[Orange,thick] (00001) -- (10002) -- (10203) -- (13204);
\draw[Orange,thick] (00001) -- (00021) -- (00132) -- (02143);
\draw[Orange,thick] (00100) -- (10200) -- (13200) -- (13204);
\draw[Orange,thick] (00100) -- (00102) -- (00132) -- (10243);
\draw[Orange,thick] (10000) -- (10200) -- (10203) -- (10243);
\draw[Orange,thick] (00021) -- (10032);
\draw[Orange,thick] (00102) -- (10203);
\draw[Orange,thick] (02100) -- (13200);
\node[label=center:\scalebox{\font}{\textcolor{Orange}{$13254$}},rectangle,inner sep=0,minimum width=\bwid pt,minimum height=\bhei pt,fill=white] at (13254){};
\node[label=center:\scalebox{\font}{\textcolor{Orange}{$-3254$}},rectangle,inner sep=0,minimum width=\bwid pt,minimum height=\bhei pt,fill=white] at (02143){};
\node[label=center:\scalebox{\font}{\textcolor{Orange}{$1-254$}},rectangle,inner sep=0,minimum width=\bwid pt,minimum height=\bhei pt,fill=white] at (10243){};
\node[label=center:\scalebox{\font}{\textcolor{Orange}{$132-4$}},rectangle,inner sep=0,minimum width=\bwid pt,minimum height=\bhei pt,fill=white] at (13204){};
\node[label=center:\scalebox{\font}{\textcolor{Orange}{$--254$}},rectangle,inner sep=0,minimum width=\bwid pt,minimum height=\bhei pt,fill=white] at (00132){};
\node[label=center:\scalebox{\font}{\textcolor{Orange}{$-32-3$}},rectangle,inner sep=0,minimum width=\bwid pt,minimum height=\bhei pt,fill=white] at (02103){};
\node[label=center:\scalebox{\font}{\textcolor{Orange}{$1--54$}},rectangle,inner sep=0,minimum width=\bwid pt,minimum height=\bhei pt,fill=white] at (10032){};
\node[label=center:\scalebox{\font}{\textcolor{Orange}{$1-2-4$}},rectangle,inner sep=0,minimum width=\bwid pt,minimum height=\bhei pt,fill=white] at (10203){};
\node[label=center:\scalebox{\font}{\textcolor{Orange}{$132--$}},rectangle,inner sep=0,minimum width=\bwid pt,minimum height=\bhei pt,fill=white] at (13200){};
\node[label=center:\scalebox{\font}{\textcolor{Orange}{$---54$}},rectangle,inner sep=0,minimum width=\bwid pt,minimum height=\bhei pt,fill=white] at (00021){};
\node[label=center:\scalebox{\font}{\textcolor{Orange}{$--2-4$}},rectangle,inner sep=0,minimum width=\bwid pt,minimum height=\bhei pt,fill=white] at (00102){};
\node[label=center:\scalebox{\font}{\textcolor{Orange}{$-32--$}},rectangle,inner sep=0,minimum width=\bwid pt,minimum height=\bhei pt,fill=white] at (02100){};
\node[label=center:\scalebox{\font}{\textcolor{Orange}{$1---4$}},rectangle,inner sep=0,minimum width=\bwid pt,minimum height=\bhei pt,fill=white] at (10002){};
\node[label=center:\scalebox{\font}{\textcolor{Orange}{$1-2--$}},rectangle,inner sep=0,minimum width=\bwid pt,minimum height=\bhei pt,fill=white] at (10200){};
\node[label=center:\scalebox{\font}{\textcolor{Orange}{$----4$}},rectangle,inner sep=0,minimum width=\bwid pt,minimum height=\bhei pt,fill=white] at (00001){};
\node[label=center:\scalebox{\font}{\textcolor{Orange}{$--2--$}},rectangle,inner sep=0,minimum width=\bwid pt,minimum height=\bhei pt,fill=white] at (00100){};
\node[label=center:\scalebox{\font}{\textcolor{Orange}{$1----$}},rectangle,inner sep=0,minimum width=\bwid pt,minimum height=\bhei pt,fill=white] at (10000){};
\node[label=center:\scalebox{\font}{\textcolor{Orange}{$\hat{0}$}},rectangle,inner sep=0,minimum width=\bwid pt,minimum height=\bhei pt,fill=white] at (0){};
\end{tikzpicture}
\caption{The interval $[1,13254]$ (top) of the permutation poset, the poset $\hat{R}(1,13254)$ (middle) and the fibres $\hat{R}^*(132,13254),\hat{R}^*(12,13254),\hat{R}^*(21,13254)$ and $\hat{R}^*(1,13254)$ (bottom).}\label{fig:big}
\end{figure}

\begin{ex}
In Figure~\ref{fig:big} we can see the poset fibration applied to the interval $[1,13254]$ of the permutation poset. So we can use Theorem~\ref{thm:muEta2} to compute the M\"obius function in the following way:
\begin{align*}
\mu(1,13254)&=\NE(1,13254)+\mu(1,132)\hat{\mu}(R^*(132,13254))\\&\,\,\,\,\,\,\,\,+\mu(1,12)\hat{\mu}(R^*(12,13254))+\mu(1,21)\hat{\mu}(R^*(21,13254))\\&=0+1+0+0
\end{align*}
\end{ex}

\cref{thm:muEta2} shows that the M\"obius function on closed pattern posets is intrinsically linked to the number of normal embeddings. This helps to explain why normal embeddings appear in many of the results on different pattern posets that have been studied independently.

In small examples such as \cref{fig:big} it is often simpler to compute the M\"obius function in the traditional way. However, \cref{cor:mobGen,thm:muEta2} are useful when considering large general intervals of pattern posets as the total spaces tends to behave in a more structured way that allows for easier analysis. Some examples of this are presented in \cref{sec:app}.

\subsection{Disconnected Intervals of a Pattern Poset}\label{sec:discon}
In this subsection we study the property of disconnectedness in pattern posets. Proposition 5.3 of~\cite{McSt13} gives a characterisation of when an interval of the classical permutation poset is disconnected, based on whether the set of embeddings can be split in a certain way, and we generalise this result to fully-closed pattern posets. First note that in a fully-closed pattern poset an embedding is uniquely determined by its zero set.

\begin{defn}\label{defn:ZS}
An interval $[\sigma,\pi]$ of a pattern poset, with $\rk(\sigma,\pi)\ge 2$, is \emph{zero split} (resp. \emph{rep-zero split}) if the embedding set (resp. representative embedding set) can be split into two disjoint non-empty sets $E_1$ and $E_2$ such that $Z(E_1)\cap Z(E_2)=\emptyset$, where $Z(E_i)$ is the union of the zero sets of the elements of $E_i$. We call $E_1$ and $E_2$ a \emph{zero split partition} of the embedding set. We say an interval of rank $k\le1$ is never zero split.

We say that an interval $[\sigma,\pi]$, with $\rk(\sigma,\pi)\ge2$, is \emph{strongly zero split} if there exists a zero split partition $E_1$ and $E_2$ of~$E^{\sigma,\pi}$ which satisfies the following condition: For all $\eta_1\in E_1$ and $\eta_2\in E_2$ there does not exist a pair~$z_1\in Z(\eta_1)$ and $z_2\in Z(\eta_2)$ such that the embeddings in $\pi$ with zero sets~$Z(\eta_1)\setminus\{z_1\}$ and $Z(\eta_2)\setminus\{z_2\}$ are embeddings of the same element $\lambda$ in $\pi$.
\end{defn}
\begin{ex}
Consider the interval $[41253,41627385]$ of the classical permutation poset. The embeddings are $\eta_1=41-273--$ and $\eta_2=--62-385$. If we partition the embeddings into the sets $E_1=\{\eta_1\}$ and $E_2=\{\eta_2\}$, then this is a zero split partition but not a strongly zero split partition. To see this is not a strongly zero split partition, note that $Z(\eta_1)=\{3,7,8\}$ and $Z(\eta_2)=\{1,2,5\}$. The embeddings with zero sets $Z(\eta_1)\setminus\{3\}$ and $Z(\eta_2)\setminus\{5\}$ are $416273--$ and~$--627385$, respectively, which are both embeddings of~$415263$. Therefore, the condition for a strongly zero split partition is violated.
\end{ex}

\begin{rem}
Any partition of the embeddings of a rank $1$ interval is a zero split partition, however we assume rank $1$ intervals to be non-zero split. The reason for this is that we are interested in zero split partitions because they imply disconnectivity, however a rank $1$ interval has an empty interior so cannot be disconnected.
\end{rem}

Next we give some properties of being zero split in relation to the embedding posets defined in \cref{sec:posfib}.  In an interval $[\sigma,\pi]$ of a fully-closed pattern poset the join of any two embeddings~$\alpha,\beta\in A^*(\sigma,\pi)$ is given by the embedding with the zero set~$Z(\alpha)\cap Z(\beta)$. We can use this to show that an interval $[\sigma,\pi]$ being zero split is intrinsically related to the connectedness of the posets $A(\sigma,\pi)$,~$A^*(\sigma,\pi)$ and~$[\sigma,\pi]$. Note that given any embedding $\eta\in E^{\sigma,\pi}$ there is a unique representative embedding $\rp(\eta)\in\hat{E}^{\sigma,\pi}$ obtained by moving all empty position to the left and full positions to the right in each adjacency.

\begin{lem}\label{lem:discon}
Consider an interval $[\sigma,\pi]$ of a fully-closed pattern poset, with $\rk(\sigma,\pi)\ge 2$, then the following conditions are equivalent:
\begin{enumerate}[label=(\arabic*)]
\item $[\sigma,\pi]$ is zero split,\label{lem:discon:cond1}
\item $A^*(\sigma,\pi)$ is disconnected,\label{lem:discon:cond2}
\item $[\sigma,\pi]$ is rep-zero split,\label{lem:discon:cond6}
\item $R^*(\sigma,\pi)$ is disconnected,\label{lem:discon:cond4}
\end{enumerate}
Furthermore, if $\rk(\sigma,\pi)\ge3$, then the above conditions are equivalent to:
\begin{enumerate}[label=(\arabic*),resume]
\item $A(\sigma,\pi)$ is disconnected.\label{lem:discon:cond3}
\item $R(\sigma,\pi)$ is disconnected.\label{lem:discon:cond5}
\end{enumerate}

\begin{proof}
Case \ref{lem:discon:cond1} $\implies$ \ref{lem:discon:cond2}. Suppose that $[\sigma,\pi]$ is zero split with the partition~$E_1$ and $E_2$ of $E^{\sigma,\pi}$. Let~$P_1$ and $P_2$ be the elements of $A^*(\sigma,\pi)$ that contain an element of $E_1$ and $E_2$, respectively. Note that any two atoms~$\eta_1\in E_1$ and~$\eta_2\in E_2$ have $Z(\eta_1)\cap Z(\eta_2)=\emptyset$, so their join is $\hat{1}$. Therefore,  $P_1$ and~$P_2$ are disconnected components of $A^*(\sigma,\pi)$.

Case \ref{lem:discon:cond2} $\implies$ \ref{lem:discon:cond1}. Suppose $A^*(\sigma,\pi)$ is disconnected with components~$P_1$ and~$P_2$, which have atoms~$E_1$ and~$E_2$, respectively. The join of any elements~${\eta_1\in E_1}$ and $\eta_2\in E_2$  equals~$\pi$ which implies that the intersection of their zero set is empty. Moreover, because this is true for any pair, it implies~$E_1$ and $E_2$ form a zero split partition of $E^{\sigma,\pi}$.

Case \ref{lem:discon:cond1} $\implies$ \ref{lem:discon:cond6}: Suppose $[\sigma,\pi]$ is zero split with the zero split partition $E_1$ and $E_2$. Let $r(E_i)$ be obtained by removing the non-representative embeddings from $E_i$. We claim the sets $r(E_1)$ and $r(E_2)$ form a rep-zero split partition of the set of representative embeddings. We know that the intersection of zero sets is empty because $Z(r(E_1))\cap Z(r(E_2))\subseteq Z(E_1)\cap Z(E_2)=\emptyset$. However, we must check that $r(E_1)$ and $r(E_2)$ are non-empty. We can get from $\eta$ to $\rp(\eta)$ by a sequence $\eta=\alpha_1,\alpha_2,\dots,\alpha_k=\rp(\eta)$, where we can get from $\alpha_i$ to $\alpha_{i+1}$ by swapping a full position with an empty position. Therefore, as $|\pi|-|\sigma|\ge2$ we know there is at least two empty positions in each $\alpha_i$ which implies that $Z(\alpha_i)\cap Z(\alpha_{i+1})\not=\emptyset$, so $\alpha_i$ and $\alpha_{i+1}$ must be in the same part of the partition, which implies that $\rp(\eta)$ is in the same part as $\eta$. Therefore, $r(E_1)$ and $r(E_2)$ must be non-empty, so $r(E_1)$ and $r(E_2)$ is a valid rep-zero split partition of~$\hat{E}^{\sigma,\pi}$.

Case \ref{lem:discon:cond6} $\implies$ \ref{lem:discon:cond1}: Suppose $[\sigma,\pi]$ is rep-zero split, with rep-zero split partition $E_1$ and $E_2$. Let~${B_i=\{\eta\in E^{\sigma,\pi}\,|\,\rp(\eta)\in E_i\}}$, for $i=1,2$. Suppose there exists a pair ${\eta\in B_1}$ and ${\phi\in B_2}$ such that~${Z(\eta)\cap Z(\phi)\not=\emptyset}$, this implies ${Z(\rp(\eta))\cap Z(\rp(\phi))\not=\emptyset}$. However, because ${\rp(\eta)\in E_1}$ and $\rp(\phi)\in E_2$ this contradicts  $E_1$ and $E_2$ being a valid rep-zero split partition. Therefore, $Z(B_1)\cap Z(B_2)=\emptyset$ so~$B_1$ and $B_2$ form a zero split partition.

Case \ref{lem:discon:cond6} $\iff$ \ref{lem:discon:cond4}: This follows by arguments analogous to \ref{lem:discon:cond1} $\implies$ \ref{lem:discon:cond2} and \ref{lem:discon:cond2} $\implies$ \ref{lem:discon:cond1}.

Cases  \ref{lem:discon:cond2} $\implies$ \ref{lem:discon:cond3} and  \ref{lem:discon:cond4} $\implies$ \ref{lem:discon:cond5}: These follow trivially because $A(\sigma,\pi)$ (resp. $R(\sigma,\pi)$) is obtained from $A^*(\sigma,\pi)$ (resp. $R^*(\sigma,\pi)$) by removing the atoms. 

Case \ref{lem:discon:cond3} $\implies$ \ref{lem:discon:cond2}. If $A(\sigma,\pi)$ is disconnected then $A^*(\sigma,\pi)$ is connected only if there is an embedding $\eta\in E^{\sigma,\pi}$ contained in elements from both components of $A(\sigma,\pi)$. However, as the interval~$[\eta,\hat{1}]$ is a Boolean lattice, with rank greater than 2, it cannot be disconnected. Therefore, no such embedding exists so  $A^*(\sigma,\pi)$ is disconnected. 

Case \ref{lem:discon:cond5} $\implies$ \ref{lem:discon:cond4}: This follows by a similar argument to that used in the case \mbox{\ref{lem:discon:cond3} $\implies$ \ref{lem:discon:cond2}},  where the only alteration is that $[\eta,\pi]$ is a product of chains, so again it cannot be disconnected when it has rank greater than~2.
\end{proof}
\end{lem}

By \cref{lem:discon}, when looking at fully-closed pattern posets we can consider either zero splitness or rep-zero splitness. For simplicity we drop the rep prefix and simply refer to zero splitness, which can be checked by looking at either the embedding set or representative embedding set. We now use \cref{lem:discon} to consider the disconnectivity of $[\sigma,\pi]$:

\begin{prop}\label{prop:discon}
Consider an interval $[\sigma,\pi]$ of a fully-closed pattern poset, where $\rk(\sigma,\pi)\ge3$. The interval $[\sigma,\pi]$ is disconnected if and only if $[\sigma,\pi]$ is strongly zero split.
\begin{proof}
Suppose that $[\sigma,\pi]$ is strongly zero split with the partition $E_1$ and~$E_2$ of~$E^{\sigma,\pi}$. Then $A(\sigma,\pi)$ is disconnected by \cref{lem:discon}. So the only way that~$[\sigma,\pi]$ is not disconnected is if there are two embeddings $\kappa_1$ and $\kappa_2$ in separate components of~$A(\sigma,\pi)$ such that $f(\kappa_1)=f(\kappa_2)$, where $f$ is the poset fibration map. First note that if $f(\kappa_1)=f(\kappa_2)$, then for any $\phi_1\le\kappa_1$ there exists a $\phi_2\le\kappa_2$ such that $f(\phi_1)=f(\phi_2)$. Therefore, we need only consider the case that $\kappa_1$ and~$\kappa_2$ are atoms.

 So suppose $\kappa_1$ and $\kappa_2$ are atoms with zero sets $Z(\kappa_i)=Z(\eta_i)\setminus\{z_i\}$, where~$\eta_i\in E_i$ and~$z_i\in Z(\eta_i)$, for $i=1,2$. However, this implies that the embeddings with zero sets~$Z(\eta_1)\setminus\{z_1\}$ and $Z(\eta_2)\setminus\{z_2\}$ are embeddings of the same element $f(\kappa_1)$ in $\pi$, which is exactly the forbidden situation in the definition of strongly zero split. Therefore, we cannot have elements from separate components mapping to the same element, so~$[\sigma,\pi]$ is disconnected. 
 
To see the other direction suppose that $[\sigma,\pi]$ is disconnected with components $P_1$ and $P_2$, and let $E_1=f^{-1}(P_1)$ and $E_2=f^{-1}(P_2)$. As $f$ is a poset fibration we know that
 $E_1$ and $E_2$ are disconnected, so $A(\sigma,\pi)$ is disconnected and \cref{lem:discon} implies $[\sigma,\pi]$ is zero split. Moreover, no two elements from different components $E_2$ and $E_2$
 maps to the same element, which implies no two atoms map to the same element. Therefore, $[\sigma,\pi]$ is strongly zero split.
\end{proof}
\end{prop}

The proof of \cref{prop:discon} allows us to see what the disconnected components of $[\sigma,\pi]$ look like when $[\sigma,\pi]$ is strongly zero split.

\begin{cor}\label{cor:discon}
If $[\sigma,\pi]$ is an interval of a fully-closed pattern poset, with $\rk(\sigma,\pi)\ge 3$, which is strongly zero split by the partition $E_1$ and $E_2$, then~$[\sigma,\pi]$ is disconnected with components $$P_i=\{\lambda\in(\sigma,\pi)\,|\,\lambda=\pi\setminus S\text{ for some } S\subset Z(E_i)\},\,\,\,\,\,\,\,for\,\,i=1,2,$$
where $\pi\setminus S$ is obtained from $\pi$ by removing the letters $\pi_i$, for all $i\in S$. 
\end{cor}

Applying \cref{prop:discon} and \cref{cor:discon} to the classical permutation poset implies Proposition~5.3 of~\cite{McSt13}. Also note that it may be possible to derive results similar to those in this section for pattern posets that are not fully-closed, which we leave as an open problem.

\subsection{Cohen-Macaulayness Preserved by a Poset Fibration}\label{sec:shellPres}
A poset is \emph{Cohen-Macaulay} if the order complex of every interval of the poset is homotopically equivalent to a wedge of top dimensional spheres. Given a poset~$P$ and an element~$p\in P$ define the induced subposet $P_{<p}=\{q\in P\,|\,q<p\}$ and similarly define~$P_{\le p}$,~$P_{>p}$ and~$P_{\ge p}$. It was first shown in \cite{Qui78} that the Cohen-Macaulay property is preserved across a poset fibration $f:P\rightarrow Q$ if the sets $f^{-1}(Q_{\ge q})$, known as the \emph{fibres}, satisfy certain conditions. The following is a variation of these results and is the dual pure form of Theorem 5.2 of \cite{Bjo05}:
\begin{prop}\label{prop:fibCMorig}
Let $P$ and $Q$ be pure posets and let $f:P\rightarrow Q$ be a poset fibration. Assume that for all $q\in Q$  there is some $p_q\in P$ such that $f^{-1}(Q_{<q})=P_{<p_q}$ and $f^{-1}(Q_{\ge q})$ is Cohen-Macaulay. If~$P$ is Cohen-Macaulay, then $Q$ is Cohen-Macaulay.
\end{prop}

We can alter the lower ideal condition of \cref{prop:fibCMorig} to get the following result:

\begin{prop}\label{prop:fibCM}
Let $P$ and $Q$ be pure posets and let $f:P\rightarrow Q$ be a poset fibration. Assume that for all $q\in Q$  there is some $p_q\in P$ such that $Q_{<q}=f(P_{<p_q})$ and $f^{-1}(Q_{\ge q})$ is Cohen-Macaulay. If~$P$ is Cohen-Macaulay, then $Q$ is Cohen-Macaulay.
\begin{proof}
Consider the posets $Q^i=\{q\in Q|\rk(q)\le i\}\cup\{p\in P|\rk(p)>i\}$ and maps $f_i:Q^{i-1}\rightarrow Q^i$ where:
$$\alpha\le_{Q^i}\beta\iff\begin{cases}\alpha\le_Q\beta,&\mbox{ and }\rk(\alpha),\rk(\beta)\le i\\\alpha\le_P\beta,&\mbox{ and }\rk(\alpha),\rk(\beta)> i\\\alpha\le_Q f(\beta),&\mbox{ and }\rk(\alpha)\le i,\rk(\beta)> i\end{cases},$$$$f_i(q)=\begin{cases}q,&\mbox{ if }\rk(q)\not=i\\f(q),&\mbox{ if } \rk(q)=i\end{cases}.$$

Note that $Q^0=P$ and $Q^{\rk(Q)}=Q$, so $Q^0$ is Cohen-Macaulay by our assumption and we proceed by an inductive argument. Assume $Q^{i-1}$ is Cohen-Macaulay, for some $i>0$, and consider~$Q^i$. We apply \cref{prop:fibCMorig} to $f_i$. Given any $q\in Q^i$, if $\rk(q)\not=i$ then $f_i^{-1}(Q^i_{\ge q})$ is $[q,\hat{1}]$ in $Q^{i-1}$, which is an interval of a Cohen-Macaulay poset and so is Cohen-Macaulay.  If $\rk(q)=i$ then~$f_i^{-1}(Q^i_{\ge q})=f^{-1}(Q_{\ge q})$ which we assumed to be Cohen-Macaulay. Furthermore, if $\rk(q)\not=i$, then $f_i^{-1}(Q^i_{<q})=Q^{i-1}_{<q}$ and if~$\rk(q)=i$, then 
$$f_i^{-1}(Q^i_{<q})=Q^i_{<q}=f(P_{<p_q})=Q^{i-1}_{<p_q}.$$
So the conditions of \cref{prop:fibCMorig} are satisfied for $f_i$ which implies $Q^i$ is Cohen-Macaulay. Therefore, by induction $Q$ is Cohen-Macaulay.
\end{proof}
\end{prop}

We can use \cref{prop:fibCM} to consider the Cohen-Macaulay property on pattern posets. First we note that the lower ideal condition of \cref{prop:fibCM} is always satisfied for pattern posets.

\begin{lem}\label{lem:lowFib}
Let $[\sigma,\pi]$ be an interval of a pattern poset along with the poset fibration $f:A(\sigma,\pi)\rightarrow(\sigma,\pi)$, we have $(\sigma,\lambda)=f(A(\sigma,\pi)_{<\ell})$, for all ${\ell\in f^{-1}(\lambda)}$.
\begin{proof}
First note that $f^{-1}(\lambda)=E^{\lambda,\pi}$, so we need to show that given any element~$\kappa\in(\sigma,\lambda)$ and~${\ell\in E^{\lambda,\pi}}$ there is an embedding of $\kappa$ in $\pi$ that is contained in $\ell$. Let~$\phi$ be an embedding of~$\kappa$ in $\lambda$ and create an embedding $\psi$ by replacing the non-empty positions of $\ell$ with $\phi$. So $\psi$ is an embedding of $\kappa$ in $\pi$ and clearly~$\kappa\le \ell$. This completes the proof.
\end{proof}
\end{lem}
\begin{cor}\label{cor:lowFib}
Let $[\sigma,\pi]$ be an interval  of a closed   pattern poset along with the poset fibration $f:R(\sigma,\pi)\rightarrow(\sigma,\pi)$, we have $(\sigma,\lambda)=f(R(\sigma,\pi)_{<\ell})$, for every~${\ell\in f^{-1}(\lambda)}$.
\end{cor}

So applying \cref{prop:fibCM} and \cref{lem:lowFib} implies the following result:

\begin{thm}\label{thm:PPCM}
Let $[\sigma,\pi]$ be an interval of a pure pattern poset such that $A^*(\lambda,\pi)$ is Cohen-Macaulay for all $\lambda\in(\sigma,\pi)$. If $A(\sigma,\pi)$ is Cohen-Macaulay, then so is~$[\sigma,\pi]$.
\end{thm}

\begin{cor}\label{cor:CECM}
Let $[\sigma,\pi]$ be an interval  of a closed   pattern poset such that $R^*(\lambda,\pi)$ is Cohen-Macaulay for all $\lambda\in(\sigma,\pi)$. If $R(\sigma,\pi)$ is Cohen-Macaulay, then so is $[\sigma,\pi]$.
\end{cor}

A Cohen-Macaulay poset cannot contain a disconnected subposet of rank greater than $2$. Therefore, Proposition \ref{prop:discon} implies the following result:

\begin{cor}
If $[\sigma,\pi]$ is an interval of a fully-closed pattern poset and contains a strongly zero split subinterval of rank greater than $2$, then $[\sigma,\pi]$ is not Cohen-Macaulay, thus not shellable.
\end{cor}

\subsection{Shellability Preserved by a Poset Fibration}\label{sec:shell}
A poset is \emph{CL-shellable} if there is an integer labelling of the edges satisfying certain conditions, we refer the reader to \cite{Wac07} for a formal definition of CL-shellability and further background. If a poset is CL-shellable, then it is shellable. It was shown in \cite{BW83} that a poset is CL-shellable if and only if it admits a recursive atom ordering, which is defined below.

Given a bounded poset $P$, a \emph{rooted interval} is an interval~$[\alpha,\beta]$ of $P$ and a chain $c$ from $\hat{0}$ to $\alpha$, and is denoted $[\alpha,\beta]_c$.
Given any chain $c$ let $c_i$ denote the element with rank $i$ in $c$, let $c_{<i}$ denote the chain of all elements of rank less than $i$ in $c$ and let $c\cdot a$ denote the chain $c$ concatenated with the element or chain $a$. Let~$\atm{P}{\alpha}$ denote the set of atoms of $[\alpha,\hat{1}]$ in $P$.

\begin{defn}\label{defn:RAO}
A bounded poset $P$ is said to admit a \emph{recursive atom ordering (RAO)} if there is an ordering $a_1,\ldots,a_t$ of $\atm{P}{\alpha}$ for every rooted interval $[\alpha,\hat{1}]_c$ that satisfies:
\begin{enumerate}[label=\textnormal{(R\arabic*)}]
\item If $r=\rk(\alpha)\ge1$, then the elements of $\F{\alpha}{c}{P}$ must appear first in the ordering, where $\F{\alpha}{c}{P}$ contains the elements of $\atm{P}{\alpha}$ which cover an element ordered before $\alpha$ in the ordering of the atoms of $[c_{r-1},\hat{1}]_{c_{<r}}$.\label{cond1}
\item For all $i<j$ if $a_i,a_j<y$, then there is a $k<j$ and an atom $z\in\atm{P}{a_j}$ such that $y\ge z >a_k$.\label{cond2}
\end{enumerate}
\end{defn}
We drop the subscripts and superscripts from $\F{\alpha}{c}{P}$ and $\atm{P}{\alpha}$ when the context is clear.

\begin{ex}
In \cref{fig:fibshell} a recursive atom ordering of $P$ is given, where the atoms of each element are ordered by their labels. To see condition \ref{cond1} consider $\atm{}{3}=\{4,6\}$. The only element of $\Omega(3)$ is $4$, because $4$ contains $2$ which is ordered before $3$.  So the condition requires that $4$ is ordered before $6$.

To see condition \ref{cond2} consider $\atm{}{2}=\{4,5,7\}$. We have $4,7\le10$ so we require an element $k$ ordered before $7$ and an element $z\in\atm{}{7}$ with $k<z<10$. So the condition is satisfied by $k=5$ and $z=9$. However, if we switched the ordering of $5$ and $7$, then it is no longer possible to find a valid pair $k<z$, hence we would no longer have an RAO.
\end{ex}
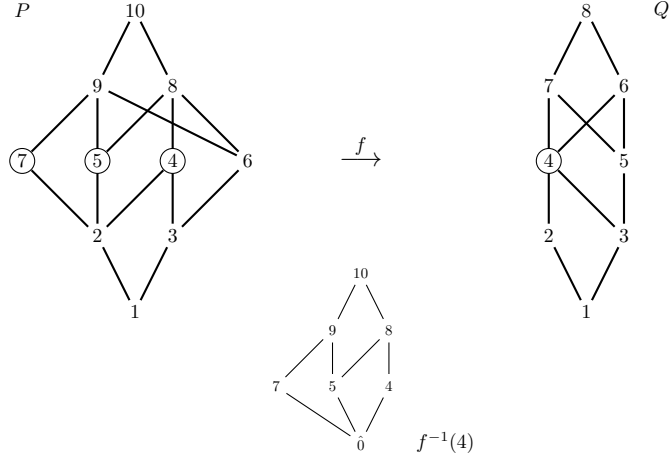
\begin{figure}[h]\centering
\begin{tikzpicture}
\def\font{.75}
\def\x{3}
\def\y{1}
\def\d{-2.5}
\node (13254) at (0+2*\x,4) {};
\node (1243) at (-.5+2*\x,3) {};
\node (2143) at (.5+2*\x,3) {};
\node (123) at (-.5+2*\x,2) {};
\node (213) at (.5+2*\x,2) {};
\node (12) at (-.5+2*\x,1) {};
\node (21) at (.5+2*\x,1) {};
\node (1) at (0+2*\x,0) {};
\draw[thick] (1) -- (12) -- (123) -- (1243) -- (13254);
\draw[thick] (1) -- (21) -- (213) -- (2143) -- (13254);
\draw[thick] (1243) -- (213);
\draw[thick] (2143) -- (123) -- (21);
\node[label=center:\scalebox{\font}{$8$},circle,fill=white] at (13254){};
\node[label=center:\scalebox{\font}{$7$},circle,fill=white] at (1243){};
\node[label=center:\scalebox{\font}{$6$},circle,fill=white] at (2143){};
\node[label=center:\scalebox{\font}{$4$},circle,fill=white,draw=black] at (123){};
\node[label=center:\scalebox{\font}{$5$},circle,fill=white] at (213){};
\node[label=center:\scalebox{\font}{$2$},circle,fill=white] at (12){};
\node[label=center:\scalebox{\font}{$3$},circle,fill=white] at (21){};
\node[label=center:\scalebox{\font}{$1$},circle,fill=white] at (1){};
\node (Q) at (1+2*\x,4){\scalebox{0.75}{$Q$}};
\node (arrow) at (0+\x,2){\scalebox{1}{$\longrightarrow$}};
\node (f) at (0+\x,2.2){\scalebox{0.75}{$f$}};
\node (10) at (0,4) {};
\node (9) at (-.5,3) {};
\node (8) at (.5,3) {};
\node (7) at (-1.5,2) {};
\node (6) at (-.5,2) {};
\node (5) at (.5,2) {};
\node (4) at (1.5,2) {};
\node (3) at (-.5,1) {};
\node (2) at (.5,1) {};
\node (1) at (0,0) {};
\draw[thick] (1) -- (2) -- (4) -- (8) -- (10) -- (9) -- (4);
\draw[thick] (1) -- (3) -- (5) -- (8);
\draw[thick] (6) -- (9) -- (7) -- (3) -- (6) -- (8);
\draw[thick] (2) -- (5);
\node[label=center:\scalebox{\font}{$10$},circle,fill=white] at (10){};
\node[label=center:\scalebox{\font}{$9$},circle,fill=white] at (9){};
\node[label=center:\scalebox{\font}{$8$},circle,fill=white] at (8){};
\node[label=center:\scalebox{\font}{$7$},circle,fill=white,draw=black] at (7){};
\node[label=center:\scalebox{\font}{$5$},circle,fill=white,draw=black] at (6){};
\node[label=center:\scalebox{\font}{$4$},circle,fill=white,draw=black] at (5){};
\node[label=center:\scalebox{\font}{$6$},circle,fill=white] at (4){};
\node[label=center:\scalebox{\font}{$2$},circle,fill=white] at (3){};
\node[label=center:\scalebox{\font}{$3$},circle,fill=white] at (2){};
\node[label=center:\scalebox{\font}{$1$},circle,fill=white] at (1){};
\node (P) at (-1.5,4){\scalebox{0.75}{$P$}};
\def\s{0.75}
\def\font{.5625}
\node (10) at (0*\s+\x,4*\s+\d) {};
\node (9) at (-.5*\s+\x,3*\s+\d) {};
\node (8) at (.5*\s+\x,3*\s+\d) {};
\node (7) at (-1.5*\s+\x,2*\s+\d) {};
\node (6) at (-.5*\s+\x,2*\s+\d) {};
\node (5) at (.5*\s+\x,2*\s+\d) {};
\node (3) at (0*\s+\x,1*\s+\d) {};
\draw (3) -- (5) -- (8) -- (10) -- (9);
\draw (6) -- (9) -- (7) -- (3) -- (6) -- (8);
\node[label=center:\scalebox{\font}{$10$},circle,fill=white] at (10){};
\node[label=center:\scalebox{\font}{$9$},circle,fill=white] at (9){};
\node[label=center:\scalebox{\font}{$8$},circle,fill=white] at (8){};
\node[label=center:\scalebox{\font}{$7$},circle,fill=white] at (7){};
\node[label=center:\scalebox{\font}{$5$},circle,fill=white] at (6){};
\node[label=center:\scalebox{\font}{$4$},circle,fill=white] at (5){};
\node[label=center:\scalebox{\font}{$\hat{0}$},circle,fill=white] at (3){};
\node (Q) at (1.5*\s+\x,1*\s+\d){\scalebox{0.75}{$f^{-1}(4)$}};
\end{tikzpicture}
\caption{A poset fibration $f:P\rightarrow Q$ which maps $\{4,5,7\}\mapsto 4$, where the linear order $\prec$ of $P$ given by the labels, induces a shelling on $P$ and the fibres of $f$, and thus $\prec_f$ induces a shelling on $Q$ by \cref{prop:fibShellrestricted}.}\label{fig:fibshell}
\end{figure}

A \emph{linear order} $\prec$ of $P$ is an ordering of all elements of $P$. We say that a linear order $\prec$ induces an RAO on $P$  if the ordering of $\atm{}{\alpha}$ by $\prec$ satisfies the RAO conditions for all rooted intervals~$[\alpha,\hat{1}]_c$ of $P$. Given a poset fibration~$f:P\rightarrow Q$ define $r_f^\prec(q)$ as the earliest element of $f^{-1}(q)$ in the linear ordering $\prec$. Define~$\prec_f$ as the linear order on $Q$ where $a\prec_f b$ if and only if $r_f^\prec(a)\prec r_f^\prec(b)$. Moreover, given a chain $c$ in $Q$ let $r_f^\prec(c)$ be the chain ${r_f^\prec(c_1)< r_f^\prec(c_2)<\cdots}$. We drop the subscripts and superscripts from $r_f^\prec(q)$ when the context is unambiguous. 

\begin{prop}\label{prop:fibShellrestricted}
Let $P$ and $Q$ be pure posets, $f:P\rightarrow Q$ be a poset fibration and $\prec$ a linear ordering of $P$. Suppose that $f^{-1}(Q_{<q})=P_{<r(q)}$ and $\prec$ induces an RAO on $f^{-1}(Q_{>q})$, for all $q\in Q$. If $\prec$ induces an RAO on $P$, then $\prec_f$ induces an RAO on $Q$, thus $Q$ is CL-shellable.
\begin{proof}
To show this we check that both conditions of an RAO are satisfied given any rooted interval~$[\alpha,\hat{1}]_c$ in $Q$. First we check Condition~\ref{cond1}. Consider~${a_i\in \F{\alpha}{c}{Q}}$ and $a_j\not\in \F{\alpha}{c}{Q}$ and let $\lambda$ be the element ordered before $\alpha$ that is covered by $a_i$. By the condition $f^{-1}(Q_{<a_i})=P_{<r(a_i)}$ we know that~$r(\lambda)\lessdot r(a_i)$ and it is straightforward to see $r(\lambda)\prec r(\alpha)$, therefore $r(a_i)\in \F{r(\alpha)}{r(c)}{P}$. Moreover, by a similar argument, $r(a_j)\not\in \F{r(\alpha)}{r(c)}{P}$. Therefore, $r(a_i)\prec r(a_j)$ so~$a_i\prec_f a_j$ and thus Condition~\ref{cond1} is satisfied.

Now we check Condition~\ref{cond2}. Consider two elements $a_i\prec_f a_j$ in $\atm{Q}{\alpha}$ and some $y>a_i,a_j$. So,~$r(a_i)$ and $r(a_j)$ are atoms of $f^{-1}(Q_{>\alpha})$ and by the condition $f^{-1}(Q_{<y})=P_{<r(y)}$ we know that~$r(y)>r(a_i),r(a_j)$. Moreover, because $\prec$ induces an RAO on $f^{-1}(Q_{>\alpha})$, there exists an atom~$\hat{a}$ of $f^{-1}(Q_{>\alpha})$ and element $z\in Q$ with $\hat{a}\prec r(a_j)$ and $\hat{a}\lessdot z\le r(y)$. Therefore, $f(\hat{a})\prec_f a_j$ and~$f(\hat{a})\lessdot f(z)\le y$, so Condition~\ref{cond2} is satisfied.
\end{proof}
\end{prop}

See \cref{fig:fibshell} for an example of \cref{prop:fibShellrestricted}. We also get the following result, whose proof we omit as it follows by an argument analogous to that used in the proof of \cref{prop:fibCM}.

\begin{prop}\label{prop:fibShellrestricted2}
Let $P$ and $Q$ be pure posets, $f:P\rightarrow Q$ a poset fibration and $\prec$ a linear ordering of $P$. Suppose that $Q_{<q}=f(P_{<r(q)})$ and $\prec$ induces an RAO in $f^{-1}(Q_{>q})$, for all $q\in Q$. If~$\prec$ induces an RAO in $P$, then $\prec_f$ induces an RAO in $Q$, thus $Q$ is CL-shellable.
\end{prop}

Applying \cref{prop:fibShellrestricted2} and \cref{lem:lowFib} to pattern posets gives the following result:
\begin{thm}\label{thm:PPshell}
Consider an interval $[\sigma,\pi]$ of a pure pattern poset $P$.
\begin{enumerate}[label=(\alph*)]
\item If $A(\sigma,\pi)$ has a linear order which induces an RAO on $A(\sigma,\pi)$ and\linebreak $A^*(\lambda,\pi)$ for all ${\lambda\in(\sigma,\pi)}$, then $[\sigma,\pi]$ is CL-shellable.\label{PPa}
\item Suppose $P$ is a closed pattern poset. If $R(\sigma,\pi)$ has a linear order which induces an RAO on~$R(\sigma,\pi)$ and $R^*(\lambda,\pi)$ for all $\lambda\in(\sigma,\pi)$, then $[\sigma,\pi]$ is CL-shellable.\label{PPb}
\end{enumerate}
\end{thm}

We believe it is possible to relax the conditions in \cref{prop:fibShellrestricted} to consider CL-shellings not induced by a linear order. So we leave the following as an open question:
\begin{que}
When can CL-shellability be preserved by a poset fibration?
\end{que}

\section{Applications}\label{sec:app}

In this section we use the results from Section~\ref{sec:results} to examine the poset of words with subword order and the consecutive permutation poset. First we introduce a lemma that proves useful. To show that~$R^*(\sigma,\pi)$ and $R(\sigma,\pi)$ are shellable using a recursive atom ordering does not require that we check every rooted interval $[\lambda,\pi]_c$. In fact it suffices to prove there is an ordering of the atoms of~$R^*(\sigma,\pi)$ which satisfies Condition~\ref{cond2}.

\begin{lem}\label{lem:satCond2}
Consider an interval $[\sigma,\pi]$ of a closed pattern poset $P$. If there is an ordering $\prec$ of~$\hat{E}^{\sigma,\pi}$ which satisfies Condition~\ref{cond2}, then $R^*(\sigma,\pi)$ and $R(\sigma,\pi)$ are shellable. Moreover, if $P$ is a fully-closed pattern poset, then
 $$\mu(\hat{R}^*(\sigma,\pi))=(-1)^{|\pi|-|\sigma|-1}|V(\sigma,\pi)|,$$
where $V(\sigma,\pi)$ is the set of embeddings $\eta\in\hat{E}^{\sigma,\pi}$ such that $\at(\eta)=\F{\eta}{}{}$.
\begin{proof}
Given any $\alpha\in R^*(\sigma,\pi)$ we refer to the \emph{filling position} of each element $\beta\in\at(\alpha)$ as the position of the letter increased to get from~$\alpha$ to~$\beta$. First we show that $R^*(\sigma,\pi)$ is shellable. By \cref{rem:prodChains} we know that~$[\alpha,\pi]$ is isomorphic to a product of chains. Therefore, every pair in~$\at(\alpha)$ are covered by their join, so any ordering of the atoms satisfies  Condition~\ref{cond2}. Define an RAO on $R^*(\sigma,\pi)$ in the following way. Consider any rooted element $\phi_c$ in~$R^*(\sigma,\pi)$. If $\rk(\phi)=0$, then order~$\at(\phi)$ according to $\prec$. If $\rk(\phi)=1$ order the elements of~$\F{\phi}{}{}$ in increasing order of the filling position and then the remaining elements in any order. If $\rk(\phi)>1$, then order $\at(\phi)$ by the order of the filling positions induced by the ordering of $\at(c_1)$. It is straightforward to see that this ordering satisfies Conditions~\ref{cond1} and~\ref{cond2}, so we have an RAO of $R^*(\sigma,\pi)$, so it is shellable. Furthermore,~$R(\sigma,\pi)$ is obtained by removing the atoms of $R^*(\sigma,\pi)$, and so $R(\sigma,\pi)$ is shellable by \cite[Theorem 8.1]{BW83}.

Now we show the M\"obius function result. Using the RAO we defined on $R^*(\sigma,\pi)$ a chain $c$ is decreasing if the empty positions of $c_1$ are filled in the reverse order of $\at(c_1)$ and $c_2$ is in $\F{c_1}{}{}$, which implies every atom of $\at(c_1)$ must be in $\F{c_1}{}{}$. So the number of decreasing chains is $|V(\sigma,\pi)|$.
\end{proof}
\end{lem}
\begin{cor}\label{cor:satCond2}
Consider an interval $[\sigma,\pi]$ of a fully-closed pattern poset. If there is an ordering of $E^{\sigma,\pi}$ which satisfies Condition~\ref{cond2}, then $A^*(\sigma,\pi)$ and $A(\sigma,\pi)$ are shellable.
\begin{proof}
By \cref{rem:closed} we know that in a fully-closed pattern poset $[\eta,\pi]$ is isomorphic to a boolean lattice, for all $\eta\in E^{\sigma,\pi}$. Therefore, the proof follows by the same argument used to prove shellability in \cref{lem:satCond2}.
\end{proof}
\end{cor}

\subsection{Poset of Words With Subword Order}
It was shown in \cite{Bjo90} that any interval $[u,w]$ of the poset of words with subword order is shellable, thus Cohen-Macaulay, and the M\"obius function equals the number of normal embeddings with sign given by the rank. In this section we give an alternative proof of Cohen-Macaulayness and the M\"obius function result on this poset.  

 Note that the poset of words with subword order is a fully-closed pattern poset and that the definition of normal embedding given by Bj\"orner is equivalent to \cref{defn:normal} when applied to this poset. Given a pair of positions $i$ and~$j$ in $\eta\in R^*(u,w)$ that are empty and non-empty, respectively, then \emph{moving} $i$ to~$j$ means setting the position $i$ as empty and the position $j$ as non-empty.

\begin{prop}
Consider an interval $[u,w]$ of the poset of words with subword order. The interval~$[u,w]$ is Cohen-Macaulay and $$\mu(u,w)=(-1)^{|w|-|u|}\NE(u,w).$$
\begin{proof}
The poset of words with subword order is a fully-closed pattern poset. So first we show that $R(u,w)$ and $R^*(u,w)$ are shellable. By \cref{lem:satCond2} it suffices to order the elements of $\hat{E}^{u,w}$ in a way satisfying Condition~\ref{cond2}. Define the \emph{position word} of an embedding as the non-empty positions listed in increasing order and the order $\prec$ on~$\hat{E}^{u,w}$ as the lexicographic order on the position words. To show that $\prec$ satisfies Condition~\ref{cond2}, consider any two embeddings~$\eta_i,\eta_j\in\hat{E}^{u,w}$, with~$\eta_i\prec\eta_j$, and some $y>\eta_i,\eta_j$. Let $a$ (resp.~$b$) be the leftmost non-empty position of $\eta_i$ (resp. $\eta_j)$ that is empty in~$\eta_j$ (resp.~$\eta_i$). The $a$'th letter of $\eta_i$ and $b$'th letter of $\eta_j$ correspond to the same letter in~$u$. Therefore, moving $b$ to $a$ in $\eta_j$ gives a valid embedding $\eta_k$ with~$\eta_k\prec\eta_j$. Moreover, let $z\in\at(\eta_j)$ be the embedding obtained by filling $a$ in~$\eta_j$, then~$\eta_k\lessdot z \le y$, so Condition~\ref{cond2} is satisfied. So $[u,w]$ is Cohen-Macaulay by \cref{cor:CECM}.

Next we consider the M\"obius function result using \cref{lem:satCond2}. Consider any embedding $\eta$  and the embedding $\phi$ obtained by filling the rightmost empty position $i$ of $\eta$, then $\phi\not\in \F{\eta}{}{}$. To see this note that if $\phi\in \F{\eta}{}{}$ then there is an element $\psi\prec\eta$  with $\psi\lessdot\phi$, where $\psi$ is obtained from~$\eta$ by moving a letter $j>i$ to~$i$. However, if $j$ is in the same adjacency as $i$, then $\psi$ would not be representative, so there must be a letter with a different value between~$i$ and $j$. However, this means we are moving~$j$ across a letter with a different value, so the order of the letters is different which implies $\psi$ is not an embedding of $u$. Therefore,~$\phi\not\in \F{\eta}{}{}$, so $\F{\eta}{}{}\not=\at(\eta)$, for any $\eta\in\hat{E}^{u,w}$, so by \cref{lem:satCond2}~$\mu(\hat{R}^*(u,w))=0$. Moreover, as this is true for any interval $[u,w]$ the M\"obius function result follows from \cref{thm:muEta2}.
\end{proof}
\end{prop}

\subsection{Consecutive Permutation Poset}
The M\"obius function and topology of the consecutive permutation poset has been studied in \cite{BFS11,SW12,EM15}, and a formula for the M\"obius function has been developed. To state this formula we first introduce some notation. A permutation is \emph{monotone} if it is of the form $12\ldots n$ or $n\ldots21$. A permutation~$\sigma$ is a \emph{prefix} of $\pi$ if $\pi_1\ldots\pi_{|\sigma|}$ is an occurrence of $\sigma$, similarly define a \emph{suffix}, and $\sigma$ is \emph{bifix} of $\pi$ if it is both a prefix and suffix of $\pi$. The \emph{exterior} of $\pi$, denoted~$x(\pi)$, is the longest bifix of~$\pi$ and the interior of $\pi$, denoted $i(\pi)$, is $\pi_2\ldots\pi_{|\pi|-1}$.

\begin{thm}\label{thm:conMob}\cite[Theorem 1.1]{SW12}
The M\"obius function of any interval $[\sigma,\pi]$ of the consecutive pattern poset is:
$$\def\b{11}\mu(\sigma,\pi)=\begin{cases}\mu(\sigma,x(\pi)),&\hskip -\b pt\mbox{ if }|\pi|-|\sigma|>2\mbox{ and }\sigma\le x(\pi)\not\le i(\pi)\\1,&\hskip -\b pt\mbox{ if }|\pi|-|\sigma|=2,\,\pi\mbox{ is non-monotone }\&\,\sigma\in\{i(\pi),x(\pi)\}\\(-1)^{|\pi|-|\sigma|},&\hskip -\b pt\mbox{ if }|\pi|-|\sigma|<2\\0,&\hskip -\b pt\mbox{ otherwise }\end{cases}.$$
\end{thm}

The consecutive pattern poset is a non-closed pattern poset. We can use \cref{cor:mobGen} to provide an alternative proof of \cref{thm:conMob}. So to compute~$\mu(\sigma,\pi)$ we need to know $\mu(\eta,\pi)$, for each $\eta\in E^{\sigma,\pi}$, and $\mu(\hat{A}^*(\lambda,\pi))$, for all $\lambda\in[\sigma,\pi)$.  Note that given any embedding $\eta$ of $\sigma$ in $\pi$ there are at most two positions that can be filled in $\eta$, the positions immediately left and right of the occurrence. So any element $\phi\ge\eta$ is obtained by a sequence of left/right fillings. We say an embedding is a prefix embedding if the non-empty positions are the initial $k$ positions, and similarly define a suffix embedding.

\begin{lem}\label{lem:con1}
Given any interval $[\sigma,\pi]$ of the consecutive permutation poset and embedding ${\eta\in E^{\sigma,\pi}}$,  we have $$\mu(\eta,\hat{1})=\begin{cases}0,&\mbox{ if } |\pi|-|\sigma|>2,\\0,&\mbox{ if } |\pi|-|\sigma|=2\mbox{ and }\eta\mbox{ is a prefix or suffix embedding,}\\(-1)^{|\pi|-|\sigma|},&\mbox{ otherwise.}\end{cases}$$
\begin{proof}
If $|\pi|-|\sigma|<2$ the result is trivial, so suppose $|\pi|-|\sigma|\ge2$. There are at most three elements of rank $2$ in~$[\eta,\pi]$, obtained from $\eta$ in the following way: $\alpha_1$ obtained by two left fillings, $\alpha_2$ obtained by two right fillings and $\alpha_3$ obtained by a left and right filling. It is straightforward to see that $\mu(\eta,\alpha_1)=\mu(\eta,\alpha_2)=0$ and $\mu(\eta,\alpha_3)=-1$. Note that if $\alpha_1$ or $\alpha_2$ equal $\hat{1}$, then $\eta$ is a prefix or suffix embedding, so the case $|\pi|-|\sigma|=2$ is complete. If $|\pi|-|\sigma|>2$, consider any element $\phi$ with rank greater than $2$. If $\phi$ does not contain $\alpha_3$ then it must be obtained by only filling left or only filling right positions, so $[\eta,\phi]$ is a chain and thus $\mu(\eta,\phi)=0$. If $\phi$ contains $\alpha_3$, then by a simple inductive argument it can be seen that $\mu(\eta,\phi)=0$, because $\alpha_3$ contains all the elements $\kappa$ with~$\mu(\eta,\kappa)\not=0$. This completes the proof.
\end{proof}
\end{lem}

\begin{lem}\label{lem:con2}
Given any interval $[\sigma,\pi]$ of the consecutive permutation poset, we have $$\mu(\hat{A}^*(\sigma,\pi))=\begin{cases}1,&\mbox{ if }\sigma=x(\pi)\mbox{ and }\sigma\not\le i(\pi)\\0,&\mbox{ otherwise}\end{cases}.$$
\begin{proof}
Given two embeddings $\eta_1$ and $\eta_2$, let $i$ and $j$ be the leftmost and rightmost positions, respectively, that are non-empty in $\eta_1$ or $\eta_2$. Then the join of~$\eta_1$ and $\eta_2$ is obtained by setting positions $i$ through $j$ as non-empty and all other positions as empty. So if $\sigma$ is not a bifix of $\pi$ then the join of the atoms of~$A^*(\sigma,\pi)$ is less than $\pi$, so $\mu(\hat{A}^*(\sigma,\pi))=0$. 

Next note that $\pi$ cannot contain a non-exterior bifix not contained in the interior. To see this let~$\alpha$ be such a bifix, then $\alpha$ is also a bifix of $x(\pi)$, which means that $\alpha$ occurs as the suffix of the prefix occurrence of $x(\pi)$, which is in the interior, giving a contradiction.

Suppose $\sigma$ is a bifix and let $\eta_1$ and $\eta_2$ be the prefix and suffix embeddings, respectively.
The join of any set of $k\ge0$ embeddings, that doesn't contain both $\eta_1$ and $\eta_2$, contributes~$(-1)^k$ to the~$\mu(\hat{A}^*(\sigma,\pi))$ and every other element contributes $0$, by the Crosscut Theorem; see~\cite[Corollary~3.9.4]{Sta97}. Therefore, \begin{equation}\label{eq:lem:con2}\mu(\hat{A}^*(\sigma,\pi))=-\sum_{S\subseteq E^{\sigma,\pi}\setminus\{\eta_1,\eta_2\}}(-1)^{|S|}+(-1)^{|S|+1}+(-1)^{|S|+1},\end{equation}
where we get three terms from considering the sets $S$, $S\cup\{\eta_1\}$ and $S\cup\{\eta_2\}$. \cref{eq:lem:con2} equals~$1$ if $E^{\sigma,\pi}=\{\eta_1,\eta_2\}$ and $0$ otherwise, which completes the proof.
\end{proof}
\end{lem}

Combining \cref{lem:con1,lem:con2} and \cref{cor:mobGen} provides an alternative proof of \cref{thm:conMob}. Moreover, Theorem 4.3 of \cite{EM15} states that an interval of the consecutive permutation poset is shellable if and only if it has no disconnected subintervals of rank greater than $2$. It is straightforward to define a shelling on $A^*(\sigma,\pi)$ if $[\sigma,\pi]$ has no disconnected subintervals. So \cref{thm:PPCM}  can be used to provide an alternative proof of the Cohen-Macaulayness of these posets.

%
%
%
%
%
%
%
%
%
%
%
%
%
%
%
%
%

\section{Future work}\label{sec:FW}
We have introduced a general definition of a pattern poset and given some results that apply to these posets. In \cref{sec:app} we applied these results to two previously studied posets and showed that our results can provide alternative proofs for existing results on these posets. There are many other pattern posets, some of which have been previously studied and many of which have not, and applying the results and techniques we have presented here could be very helpful in the study of these posets. For example, very little is known of the Dyck path poset introduced in \cite{Bac14}, which seems to have many nice properties, such as the sign of the M\"obius function being alternating. Can we apply some of the results we have introduced to learn more about this poset?

One particular pattern poset for which there are many open problems is the classical permutation poset.  The results from \cref{sec:mobGen} imply the main results of \cite{Smith15}, and the results from \cref{sec:discon} imply some of the results in \cite{McSt13}. Whether we can apply the results from \cref{sec:shellPres,sec:shell} to determine the topology of intervals of the classical permutation poset is still open. We conjecture that if an interval $[\sigma,\pi]$ of the classical permutation pattern poset does not contain any zero split subintervals, then $[\sigma,\pi]$ is shellable. If we can find an ordering of the embeddings of such intervals which satisfies Condition~\ref{cond2},  then \cref{lem:satCond2} would prove that these intervals are Cohen-Macaulay and allow us to compute the M\"obius function of these intervals. 

Our definition of normal is not equivalent to those given for the poset of words with composition order or for generalised subword order. Therefore, the results on the M\"obius function of these posets given in~\cite{SagVat06} and \cite{McnSag12} do not follow immediately from \cref{thm:muEta2}. However, it seems reasonable to hope that with some work one could apply \cref{thm:muEta2} to provide an alternative proof of these results. Moreover, such an alternative proof might provide further insight into the structure of these posets.


\section*{Acknowledgements}
I would like express my gratitude to the referees, whose comments greatly improved the paper.

\newcommand{\etalchar}[1]{$^{#1}$}
 \newcommand{\noop}[1]{}


\begin{thebibliography}{BFPW13}

\bibitem[BBF{\etalchar{+}}14]{Bac14}
Axel Bacher, Antonio Bernini, Luca Ferrari, Benjamin Gunby, Renzo Pinzani, and
  Julian West.
\newblock The {D}yck pattern poset.
\newblock {\em Discrete Mathematics}, 321:12--23, 2014.

\bibitem[BF14]{BF14}
Antonio Bernini and Luca Ferrari.
\newblock Vincular pattern posets and the {M}\" obius function of the
  quasi-consecutive pattern poset.
\newblock {\em arXiv preprint 1410.6046}, 2014.

\bibitem[BFPW13]{BLPW13}
Antonio Bernini, Luca Ferrari, Renzo Pinzani, and Julian West.
\newblock Pattern-avoiding dyck paths.
\newblock In {\em 25th International Conference on Formal Power Series and
  Algebraic Combinatorics (FPSAC 2013)}, pages 683--694. Discrete Mathematics
  and Theoretical Computer Science, 2013.

\bibitem[BFS11]{BFS11}
Antonio Bernini, Luca Ferrari, and Einar Steingr{\'i}msson.
\newblock The {M}{\"o}bius function of the consecutive pattern poset.
\newblock {\em The Electronic Journal of Combinatorics}, 18(1):P146, 2011.

\bibitem[BJJS11]{BJJS11}
Alexander Burstein, V{\'{\i}}t Jel{\'{\i}}nek, Eva Jel{\'{\i}}nkov{\'a}, and
  Einar Steingr{\'{\i}}msson.
\newblock The {M}\"obius function of separable and decomposable permutations.
\newblock {\em Journal of Combinatorial Theory. Series A}, 118(8):2346--2364,
  2011.

\bibitem[Bj{\"o}80]{Bjo80}
Anders Bj{\"o}rner.
\newblock Shellable and {C}ohen-{M}acaulay partially ordered sets.
\newblock {\em Transactions of the American Mathematical Society},
  260(1):159--183, 1980.

\bibitem[Bj{\"o}90]{Bjo90}
Anders Bj{\"o}rner.
\newblock The {M}{\"o}bius function of subword order.
\newblock {\em Institute for Mathematics and its Applications}, 19:118, 1990.

\bibitem[Bj{\"o}93]{Bjo93}
Anders Bj{\"o}rner.
\newblock The {M}{\"o}bius function of factor order.
\newblock {\em Theoretical Computer Science}, 117(1):91--98, 1993.

\bibitem[BW83]{BW83}
Anders Bj{\"o}rner and Michelle Wachs.
\newblock On lexicographically shellable posets.
\newblock {\em Transactions of the American Mathematical Society},
  277(1):323--341, 1983.

\bibitem[BWW05]{Bjo05}
Anders Bj{\"o}rner, Michelle Wachs, and Volkmar Welker.
\newblock Poset fiber theorems.
\newblock {\em Transactions of the American Mathematical Society},
  357(5):1877--1899, 2005.

\bibitem[EM15]{EM15}
Sergi Elizalde and Peter~R.W. McNamara.
\newblock The structure of the consecutive pattern poset.
\newblock {\em arXiv preprint 1508.05963}, 2015.

\bibitem[Kit11]{Kit11}
Sergey Kitaev.
\newblock {\em Patterns in Permutations and Words}.
\newblock Monographs in Theoretical Computer Science. An EATCS Series.
  Springer, Heidelberg, 2011.

\bibitem[MS12]{McnSag12}
Peter R.~W. McNamara and Bruce~E. Sagan.
\newblock The {M}{\"o}bius function of generalized subword order.
\newblock {\em Advances in Mathematics}, 229(5):2741--2766, 2012.

\bibitem[MS15]{McSt13}
Peter R.~W. McNamara and Einar Steingr{\'\i}msson.
\newblock On the topology of the permutation pattern poset.
\newblock {\em Journal of Combinatorial Theory, Series A}, 134:1--35, 2015.

\bibitem[Qui78]{Qui78}
Daniel Quillen.
\newblock Homotopy properties of the poset of nontrivial {$p$}-subgroups of a
  group.
\newblock {\em Advances in Mathematics}, 28(2):101--128, 1978.

\bibitem[Smi14]{Smith13}
Jason~P. Smith.
\newblock On the {M}\"obius function of permutations with one descent.
\newblock {\em The Electronic Journal of Combinatorics}, 21:2.11, 2014.

\bibitem[Smi16]{Smith14}
Jason~P. Smith.
\newblock Intervals of permutations with a fixed number of descents are
  shellable.
\newblock {\em Discrete Mathematics}, 339(1):118 -- 126, 2016.

\bibitem[Smi17]{Smith15}
Jason~P. Smith.
\newblock A formula for the {M}\"obius function of the permutation poset based
  on a topological decomposition.
\newblock {\em To appear in Advances in Applied Mathematics}, 2017.

\bibitem[ST10]{SteTen10}
Einar Steingr{\'{\i}}msson and Bridget~Eileen Tenner.
\newblock The {M}\"obius function of the permutation pattern poset.
\newblock {\em Journal of Combinatorics}, 1(1):39--52, 2010.

\bibitem[Sta12]{Sta97}
Richard~P. Stanley.
\newblock {\em Enumerative Combinatorics. {V}ol. 1}.
\newblock Cambridge Studies in Advanced Mathematics. Cambridge University
  Press, second edition, 2012.

\bibitem[SV06]{SagVat06}
Bruce~E. Sagan and Vincent Vatter.
\newblock The {M}\"obius function of a composition poset.
\newblock {\em Journal of Algebraic Combinatorics}, 24(2):117--136, 2006.

\bibitem[SW12]{SW12}
Bruce~E Sagan and Robert Willenbring.
\newblock Discrete {M}orse theory and the consecutive pattern poset.
\newblock {\em Journal of Algebraic Combinatorics}, 36(4):501--514, 2012.

\bibitem[Wac07]{Wac07}
Michelle~L. Wachs.
\newblock Poset topology: {T}ools and applications.
\newblock In {\em Geometric Combinatorics}, volume~13 of {\em IAS/Park City
  Math. Ser.}, pages 497--615. Amer. Math. Soc., 2007.

\bibitem[Wal81]{Wal81}
James~W Walker.
\newblock Homotopy type and {E}uler characteristic of partially ordered sets.
\newblock {\em European Journal of Combinatorics}, 2(4):373--384, 1981.

\end{thebibliography}
\end{document}